%% file: 0main.tex
\documentclass[11pt,a4paper]{article}
\usepackage{amsmath,amsthm,amsthm,mathrsfs}
\usepackage{fancyhdr}
\usepackage[british]{babel}
\usepackage{amssymb}
\usepackage{latexsym}
\usepackage{graphicx}
\usepackage{epstopdf}
\usepackage{float}
\usepackage{color}
\usepackage[natural]{xcolor}
\usepackage{textgreek}
\usepackage{float,color,fancybox,shapepar,setspace,hyperref}
\usepackage{enumitem}
\theoremstyle{plain}
\usepackage{todonotes}
\newtheorem{theorem}{Theorem}
\newtheorem{claim}[theorem]{Claim}
\newtheorem{lemma}[theorem]{Lemma}
\newtheorem{conjecture}[theorem]{Conjecture}
\numberwithin{theorem}{section}
\newtheorem{corollary}[theorem]{Corollary}
\newtheorem{proposition}[theorem]{Proposition}

\newtheorem{question}[theorem]{Question}
\newtheorem{definition}[theorem]{Definition}
\newtheorem{fact}[theorem]{Fact}

\newcounter{propcounter}
\newcommand{\eps}{\varepsilon}
\newcommand{\<}{\subseteq}

\title{\LARGE Topological cliques in sparse expanders}
\author{
Xia Wang\thanks{School of Mathematics, Shandong University, Jinan, China. Email: {\tt xiawang@mail.sdu.edu.cn}.}
\and
Donglei Yang\thanks{School of Mathematics, Shandong University, Jinan, China. Email: {\tt dlyang@sdu.edu.cn}. Supported by Natural Science Foundation of China (12101365) and Natural Science Foundation of Shandong Province (ZR2021QA029).}
\and
Fan Yang\thanks{Data Science Institute, Shandong University, Jinan, China. Email: {\tt fyang@sdu.edu.cn}. Supported by Natural Science Foundation of China (12301447), Natural Science Foundation of Shandong Province (ZR2024QA056) and the China Postdoctoral Science Foundation (12570073310023).}
\and
Haotian Yang\thanks{School of Mathematics, Shandong University, Jinan, China. Email: {\tt 202017091012@mail.sdu.edu.cn.}
}
}
\date{}

\begin{document}
\maketitle

\input{1abctrac}
\input{2introduction}
\input{3preliminary}
\input{4.1sparse}	
\input{4.2dense}

\input{5balancedimmersion}
\input{6subdivision}

\input{7.0Acknowledgements}
\bibliographystyle{IEEEtranS}
\bibliography{7reference}
\input{8appendix}

\end{document}

%% file: 1abctrac.tex
\begin{abstract}
In the paper, we focus on embedding clique immersions and subdivisions within sparse expanders, and we derive the following main results:
 \begin{itemize}
 \item For any $0< \eta< 1/2$, there exists $K>0$ such that for sufficiently large $n$, every $(n,d,\lambda)$-graph $G$ contains a $K_{(1-5\eta)d}$-immersion when $d\geq K\lambda$.
  \item For any $\varepsilon>0$ and $0<\eta <1/2$, the following holds for sufficiently large $n$. Every $(n,d,\lambda)$-graph $G$ with $2048\lambda/\eta^2<d\leq \eta n^{1/2-\varepsilon}$ contains a $K_{(1-\eta)d}^{(\ell)}$-subdivision, where $\ell = 2 \left\lceil \log(\eta^2n/4096)\right\rceil + 5$.
 \item There exists $c>0$ such that the following holds for sufficiently large $d$. If $G$ is an $n$-vertex graph with average degree $d(G)\geq d$, then $G$ contains a $K_{c d}^{(\ell)}$-immersion for some $\ell\in \mathbb{N}$.
\end{itemize}

In 2018, Dvo{\v{r}}{\'a}k and Yepremyan asked whether every graph $G$ with $\delta(G)\geq t$ contains a $K_t$-immersion. Our first result shows that it is asymptotically true for $(n,d,\lambda)$-graphs when $\lambda=o(d)$. In addition, our second result extends a result of  Dragani{\'c}, Krivelevich and Nenadov on balanced subdivisions. The last result generalises a result of DeVos, Dvo{\v{r}}{\'a}k, Fox, McDonald, Mohar, Scheide on $1$-immersions of large cliques in dense graphs.

\end{abstract}


%% file: 2introduction.tex
\section{Introduction}
\par A graph $G$ contains an \emph{H-minor} if there are $|V(H)|$ disjoint subsets $\{T_v\}_{v\in V(H)}$ of $V(G)$ such that $T_v$ induces a connected graph and there is an edge between $T_u$ and $T_v$ in $G$ for every $uv\in E(H)$.
In 1943, Hadwiger \cite{hadwiger1943klassifikation} proposed a conjecture that every graph $G$ with $\chi(G)\geq t$ contains a $K_t$-minor, and proved the cases $t\leq 4$. Additionally, the cases $t\in\{5,6\}$ had been shown to be equivalent to the Four-Color Theorem, which were respectively proved by Wagner \cite{wagner1937eigenschaft} and Robertson, Seymour and Thomas \cite{robertson1993hadwiger}, while the cases $t\geq7 $ are still open.

A graph $G$ contains an \emph{H-subdivision} if there is an injective mapping $\phi: V(H)\to V(G)$ such that for every edge $uv\in E(H)$, there is a path $P_{uv}$ in $G$ connecting vertices $\phi(u)$ and $\phi(v)$, and all paths $P_{uv}$ in $G$ ($uv\in E(H)$) are pairwise internally vertex disjoint. 
It is easy to see that if $G$ contains an $H$-subdivision, then $G$ contains an $H$-minor. In 1940's, Haj\'os \cite{hajos1961uber} proposed a stronger conjecture which says that every graph $G$ with $\chi(G)\geq t$ contains a $K_{t}$-subdivision.
Dirac \cite{Dirac1952APO} proved the cases $t\leq 4$. But the cases $t\geq 7$ were disproved by Catlin \cite{catlin1978bound}. The cases $t\in \{5,6\}$ are still open.

\par Another related notion of immersion was first introduced by Nash-Williams \cite{nash1965well} in 1965. 
A graph $G$ contains an \emph{$H$-immersion} if there is an injective mapping $\phi: V(H)\to V(G)$ such that for every edge $uv\in E(H)$, there is a path $P_{uv}$ in $G$ connecting $\phi(u)$ and $\phi(v)$, and all paths $P_{uv} \in G$ ($uv\in E(H)$) are pairwise edge disjoint. The vertices in $\{\phi(v) \mid v\in V(H)\}$ are \emph{branch vertices}.
In particular, if $G$ contains an $H$-immersion such that all branch vertices are not the internal vertices of paths $P_{uv}$ for every $uv\in E(H)$, then $G$ contains a \emph{strong $H$-immersion}.

By the above definitions, we can see that if $G$ contains an $H$-subdivision, then it also contains an $H$-immersion. Robertson and Seymour proved that graphs are well quasi-ordered by the immersion relation \cite{robertson2010graph} based on their minors project \cite{robertson2004graph}. 
Analogous to Haj\'os' conjecture, Lescure and Meyniel~\cite{lescure41problem}, and independently, Abu-Khzam and Langston~\cite{abu2003graph} proposed the following conjecture.

\begin{conjecture}{\rm\cite{lescure41problem, abu2003graph}}\label{con1}
	Every graph $G$ with $\chi(G)\geq t$ contains a $K_t$-immersion. 
\end{conjecture}

For Conjecture \ref{con1}, the cases $t\leq 4$ follow directly from the result of Dirac \cite{Dirac1952APO}. DeVos, Kawaraba-yashi, Mohar and Okamura \cite{devoss2010immersing} solved the cases $5\leq t\leq 7$. The cases $t\geq 8$ are still open. 

\subsection{Extremal density for topological cliques}
\par The above conjectures are all asking for large topological cliques under the chromatic number condition. We know that a graph with chromatic number at least $r$ must contains a subgraph with minimum degree at least $r-1$, so many scholars turned to study the above problems under the degree condition, and obtained many classical results. For minors, Kostochka \cite{kostochka1984lower}, and independently, Thomason \cite{thomason2001extremal} proved that a graph $G$ with average degree $\Omega(d\sqrt{\log d})$ contains a $K_d$-minor, and this bound is best possible. 
For subdivisions, Mader initiated a foundamental problem of determining the largest clique subdivisions in graphs with given average degree. Bollob\'as and Thomason \cite{bollobas1998proof} and independently, Koml\'os and Szemer\'edi \cite{komlos1996topological} proved that there is a constant $C>0$ such that every graph with average degree at least $Cd$ contains a $K_{\sqrt{d}}$-subdivision, and this is tight up to the constant factor $C$. Similarly, for immersions, DeVos et al. \cite{devoss2010immersing} proposed a question as follows:  let $t$ be a positive integer, and can we find the smallest value $f(t)$ such that every graph $G$ with $\delta(G)\geq f(t)$ contains a $K_t$-immersion? 
For this question, a trivial lower bound is $f(t)\geq t-1$. Lescure and Meyniel \cite{lescure41problem} and DeVos et al. \cite{devoss2010immersing} proved that $f(t)=t-1$ for $t\leq 7$. As to $t\geq 8$, there is an infinite number of constructions \cite{devoss2010immersing,lescure41problem,collins2014constructing,devos2014minimum} indicating that $f(t)\geq t$. As for the upper bound, DeVos et al. \cite{devos2014minimum} proved that $f(t)\leq 200t$. Later, Dvo{\v{r}}{\'a}k and\ Yepremyan \cite{dvovrak2018complete} improved this bound to $f(t)\leq 11t+7$ and raised the following question.
\begin{question}{\rm\cite{dvovrak2018complete}}\label{minimumimmersion}
    Does every graph with minimum degree $t$ contains an $K_t$-immersion?
\end{question}
\noindent Up to now, the best known upper bound is $f(t)\leq 7t+7$, given by Gauthier, Le and Wollan \cite{gauthier2019forcing}.
However, Liu, Wang and Yang \cite{liu2022clique} confirmed that  Conjecture \ref{con1} and Question \ref{minimumimmersion} are asymptotically true for graphs without any fixed complete bipartite graph.

\subsection{Balanced subdivisions and immersions}
\par
In this paper, we consider topological cliques with certain length constraints. To begin with, an \emph{$H^{(\ell)}$-immersion} which is also called a \emph{balanced immersion of $H$}, is an $H$-immersion satisfying that the length of each path $P_{uv}$ ($uv\in E(H)$) is exactly $\ell+1$.
Some results of embedding immersions in graphs have been presented earlier, so what is the maximum value of $t$ if we want to embed a balanced immersion of $K_t$ to a graph $G$ with average degree $d(G)=d$? For this question, DeVos et al.~\cite{devos2014minimum} proved that for dense graphs on $n$ vertices with at least $2cn^2$ edges, there is a $K_{c^2 n}^{(1)}$-immersion. Our first work is to generalize this result to graphs of arbitrarily average degrees and obtain the following result.

\begin{theorem}\label{thm1}
There exists $c>0$ such that for any $d>0$, if $G$ is a graph with average degree $d(G)\geq d$, then $G$ contains a $K_{cd}^{(\ell)}$-immersion for some $\ell\in \mathbb{N}$.
\end{theorem}
Our second result concerns clique immersions in \emph{$(n,d,\lambda)$-graphs} with large spectral gap, which is known as ``spectral expanders" introduced by Alon \cite{alon1986eigenvalues}.
\begin{definition}{\rm\cite{alon1986eigenvalues}}
    An $(n,d,\lambda)$-graph is an $n$-vertex $d$-regular graph with all but the largest eigenvalues of its adjacency matrix being at most $\lambda$ in absolute value.
\end{definition}
\noindent We get the following result which shows that Conjecture \ref{con1} and Question \ref{minimumimmersion} are asymptotically true for $(n,d,\lambda)$-graphs with $\lambda=o(d)$.

\begin{theorem}\label{ndlthm}
For any $0< \eta< 1/2$, there exists $K>0$ such that for sufficiently large $d$, every $(n,d,\lambda)$-graph $G$ with $d\geq K\lambda$ contains a $K_{(1-5\eta)d}$-immersion.	
\end{theorem}


\par Similarly for a positive integer $\ell$, an \emph{$H^{(\ell)}$-subdivision} which is also called a \emph{balanced subdivision} of $H$, is an $H$-subdivision satisfying that the length of each path $P_{uv}$ $(uv\in E(H))$ is exactly $\ell + 1$. 
Motivated by Mader’s
conjecture, Thomassen proposed the following well-known conjecture \cite{thomassen1984subdivisions,thomassen1989configurations}.
\begin{conjecture}{\rm\cite{thomassen1984subdivisions,thomassen1989configurations}}\label{conthoma}
 For any constant $t\in\mathbb{N}$, there exists a $g(t)$ such that every graph $G$ with $d(G)\geq g(t)$ contains a balanced $K_t$-subdivision.   
\end{conjecture}

In 2020, Liu and Montgomery \cite{liu2020solution} confirmed Conjecture \ref{conthoma}. After that, Wang \cite{wang2021balanced} improved it to $g(t)=t^{2+o(1)}$.
Finally, Luan, Tang, Wang and Yang \cite{Luan2022BalancedSO} and Fern\'andez, Hyde, Liu, Pikhurko and Wu \cite{FERNANDEZ2023417} independently proved that $g(t)=\Theta(t^{2})$ in Conjecture \ref{conthoma}. Moreover, in \cite{Luan2022BalancedSO}, they showed that there exists an absolute constant $c>0$ such that every $C_4$-free graph with average degree at least $d$ contains a balanced $K_{cd}$-subdivision. 

As to the spectral expanders,  Dragani\'c, Krivelevich and Nenadov \cite{draganic2022rolling} proved a much better bound on the order of balanced clique subdivision as follows.
\begin{theorem}{\rm\cite{draganic2022rolling}}\label{d51}
Let $G$ be a $(n, d, \lambda)$-graph with $240\lambda < d \leq n^{1/5}/2$, and let $d_0 \geq 3$.
Then $G$ contains a balanced $K_t$-subdivision for $t = \lfloor d - 80\lambda\sqrt{d_0}\rfloor$, all the subdivided paths being of equal length $\ell$, where $\ell=O(\frac{\log n}{\log d_0})$.
\end{theorem}
Erd\H{o}s and Fajtlowitz \cite{erdHos1981conjecture} observed that in general we have $d=O(\sqrt{n})$ whenever one want to guarantee a $K_{d-o(d)}$-subdivision (space barrier). We extend Theorem \ref{d51} to nearly optimal bound on $d$ and obtain the following result. 
\begin{theorem}\label{d21}
For any $\varepsilon>0$ and $0<\eta <1/2$, the following holds for sufficiently large $n$. Every $(n,d,\lambda)$-graph $G$ with $2048\lambda/\eta^2<d\leq \eta n^{1/2-\varepsilon}$ contains a $K_{(1-\eta)d}^{(\ell)}$-subdivision, where $\ell = 2 \left\lceil \log(\eta^2n/4096)\right\rceil + 5$.
\end{theorem}

The rest of this paper will be organized as follows. In Section \ref{preliminary}, we introduce some notions and tools. In Section \ref{immersion}, we give a proof of Theorem \ref{ndlthm}. Section \ref{balancedimmersion} shows a proof of Theorem \ref{thm1} which is divided into three cases---dense, sparse, and intermediate case---and the proof of the intermediate case will be given in Subsection \ref{hpgra}. Finally, we complete the proof of Theorem \ref{d21} in Section \ref{subdivision}.

%% file: 3preliminary.tex
\section{Preliminary}\label{preliminary}
\par In subsection \ref{notation}, we give some definitions and notions. In subsection \ref{tools1} and \ref{tools2}, we shall present several essential tools about sublinear robust expander and $(n,d,\lambda)$-graphs.

\subsection{Notation}\label{notation}
\par Let $G$ be a graph with vertex set $V(G)$ and edge set $E(G)$. 
By $[V(G)]^2$ we denote the set of all 2-element subsets of $V(G)$.
The \emph{complement} $\overline{G}$ of $G$ is the graph on $V(G)$ with edge set $[V(G)]^2\setminus E(G)$. 
For a subset $X\subseteq V(G)$, $G[X]$ denotes the \emph{induced subgraph} of $G$ on $X$.
For any $X,Y\subseteq V(G)$, let $G[X,Y]$ be the induced bipartite graph in $G$ with two parts $X$ and $Y$, and let $e(X,Y)$ be the number of edges in $G[X,Y]$.
Let $G-X$ be the graph $G[V\setminus X]$. For a edge set $Y\subseteq E(G)$, let $G\setminus Y$ be the spanning graph $(V,E\setminus Y)$.
Denote the \emph{external\ neighborhood} of $X$ by $N_G (X):=\{u\in V(G)\setminus X: uv\in E(G) \text{ for some } v\in X\}$.
For two sets $A,B\subseteq V(G)$, denote the \emph{common\ neighbourhood} of $A$ and $B$ by $N_G(A,B):=N_G(A)\cap N_G(B)$. As before, we write $d(u,B)$ rather than $d(\{u\},B)$, etc.
The \emph{codegree} $d(u,v)$ of two vertices $u,v\in V(G)$ is the number of vertices in $N_G(u,v)$.
The \emph{length} of a path (or cycle) denotes the number of edges in the path (or cycle). 
We call $P=x_0 x_1\dots x_k$ an $(A,B)$-path if $V(P)\cap A=\{x_0\}$ and $V(P)\cap B=\{x_k\}$. As before, we write an $(a,B)$-path rather than an $(\{a\},B)$-path, etc. 
The \emph{distance} between $A$ and $B$, denoted by $dist(A,B)$, is the length of a shortest $(A,B)$-path in $G$.
An $S(u)$-star is a star centered at the vertex $u$, and the size of a star is the number of the edges in the star.
Let $H=(A,B;E)$ be a bipartite graph with vertex set $A\cup B$. The $\emph{edge density}$ of $H$ is $|E|/(|A||B|)$. For convenience, we say a graph $G$ is \emph{$(s, t, \gamma)$-dense} if for any subset $U\subseteq V(G)$ with $|U|\leq s$, and any $W\subseteq E(G)$ with $|W|\leq t$, we have $d(G\setminus W-U)\geq \gamma.$

Let $[n]:=\{1,2,\dots ,n\}$. When it is not essential, we omit the floors and ceilings, and all logarithms are natural in the whole paper.
\subsection{Robust sublinear expander}\label{tools1}
\par In short, an expander is a highly connected yet sparse graph that has some excellent properties of sparse random graphs, with broad applications in fields such as theoretical computer science, group theory, and combinatorics. Recently, research on expanders has attracted much attention. Back in 1996, Koml\'os and Szemer\'edi introduced the concept of sublinear expanders in their work \cite{komlos1996topological}.

\begin{definition}{\rm\cite{komlos1996topological}}
 \rm For each $\epsilon_1 > 0$ and $k > 0$, a graph $G$ is an $(\epsilon_1, k)$-expander if $|N(X)|\geq \rho(|X|)\cdot|X|$ for all $X \subseteq V(G)$ of size $k/2 \leq|X|\leq |V(G)|/2$, where $\rho(x)$ is the function
\[\rho(x) = \rho(x, \epsilon_1, k) := \left\{
\begin{aligned}
	&0 & \text{if} \quad x < \frac{k}{5},\\
	&\frac{\epsilon_1}{\log^2(\frac{15x}{k})} & \text{if} \quad x \geq \frac{k}{5}.
\end{aligned}
\right.
\]
\end{definition}

In this paper, we adopt a stronger notion of robust sublinear expander, proposed by Haslegrave, Kim and Liu \cite{haslegrave2022extremal}.

\begin{definition}{\rm\cite{haslegrave2022extremal}}
\rm A graph $G$ is an \emph{$(\epsilon_1, k)$-robust-expander} if for all subsets $X \subseteq V (G)$ of size $k/2 \leq |X| \leq |V (G)|/2$ and any $F \subseteq E(G)$ with $|F| \leq d(G)\rho(|X|)|X|$, we have that 
\[|N_{G\setminus F} (X)| \geq \rho(|X|)|X|.\]
\end{definition}

Let $0<\epsilon_1,\epsilon_2<1$, and throughout the rest of the paper, we write \[m:=\tfrac{2}{\epsilon_1}\log ^3 \tfrac{15n}{\epsilon_2 d}.\]
Note that for sufficiently large $d$, $\rho(x)=\rho(x,\epsilon_1,\epsilon_2d)$ is decreasing in the interval $[\epsilon_2 d/2, n/2]$, so we have that for every $\epsilon_2 d/2\leq x\leq n/2$, 
\begin{equation}\label{1m}
\epsilon_1>\tfrac{\epsilon_1}{\log^2(\tfrac{15}{2})}=\rho(\tfrac{\epsilon_2d}{2})\geq\rho(x)\geq \rho(n)\geq \tfrac{1}{m}.
\end{equation}
When $n/d$ is sufficiently large, we get 
\begin{equation}\label{m200}
\frac{n}{d}\geq m^{200}.
\end{equation}
If $d\geq \log^{200} n$, then 
\begin{equation}\label{m60}
d\geq m^{60}.
\end{equation}

Haslegrave et al. \cite{haslegrave2022extremal} proved that every graph $G$ contains a robust sublinear expander, which is almost as dense as $G$.

\begin{lemma}{\rm\cite{haslegrave2022extremal}}\label{re}
Let $C > 30$, $0 < \epsilon_1 \leq 1/(10C)$, $0 < \epsilon_2 < 1/2$, $d > 0$ and $\eta = C\epsilon_1/ \log3$. Then every graph $G$ with $d(G) = d$ has a subgraph $G'$ which is an $(\epsilon_1, \epsilon_2d)$-robust-expander with $d(G') \geq (1-\eta)d$ and $\delta(G') \geq d(G')/2$.
\end{lemma}

Note that every graph $G$ contains a bipartite subgraph $H$ with $d(H)\geq d(G)/2$. Then by Lemma \ref{re}, we immediately obtain the following corollary.

\begin{corollary}\label{col}
There exist $0<\epsilon_1,\epsilon_2<1$ such that every graph $G$ with $d(G)\geq 8d$ has a bipartite subgraph $G'$ which is an $(\epsilon_1, \epsilon_2d)$-robust-expander with $d(G') \geq 2d$ and $\delta(G') \geq d$.
\end{corollary}

The following lemma is essential to find a short path connecting two large sets while avoiding some vertices and edges.

\begin{lemma}[Robust small diameter,~\cite{haslegrave2022extremal}]\label{rbsdl}
Let $0 < \epsilon_1, \epsilon_2 < 1$ and $G$ be an $n$-vertex $(\epsilon_1, \epsilon_2 d)$-robust-expander. Given two sets $X_1, X_2 \subseteq V (G)$ of order $x\geq \epsilon_2 d/2$, let $U$ be a vertex set of order at most $\rho(x) x/4$ and $W$ be a subgraph with at most $d(G)\rho(x)x$ edges. Then there is an $(X_1, X_2)$-path of length at most $(2/\epsilon_1)\log ^3 (15n/(\epsilon_2 d))$ in $(G\setminus W)-U$ .
\end{lemma}

\subsection{\texorpdfstring{$(n,d,\lambda)$}{(n,d,\textlambda)}-graphs}\label{tools2}
\par The concept of $(n,d,\lambda)$-graphs is particularly significant in the study of expanders from various aspects. The edges in an $(n,d,\lambda)$-graphs are uniformly distributed, exhibiting excellent connectivity, and relevant research in this area can be found in recent
comprehensive surveys~\cite{Krivelevich2006, hoory2006expander}.

\begin{fact}\label{ndlcompletement}{\rm\cite{cvetkovic2009introduction}}
If $G$ is an $(n,d,\lambda)$-graph with $\lambda\geq \lambda_2\geq\cdots\geq\lambda_n$, then $\overline{G}$ is an $(n,n-1-d,-(\lambda_n+1))$-graph.
\end{fact}

\begin{lemma}{\rm\cite{alon2016probabilistic}}\label{Alon}
For every  $(n,d,\lambda)$-graph $G$ and for every partition of the set of vertices $V(G)$ into two disjoint subsets $B$ and $C$, it holds 
\[ e(B,C)\geq \frac{(d-\lambda)|B||C|}{n}. \]
\end{lemma}

This implies that if the second largest eigenvalue of $G$ is far from the first, then $G$ has good expansion properties. The converse of this is also true, and you can see more details in \cite{alon2016probabilistic}. Another famous one is expander mixing lemma, which intuitively states that the edges of certain $d$-regular graphs are evenly distributed throughout the graph. In particular, the number of edges between two vertex subsets $U$ and $V$ is always close to the expected number of edges between them in a random $d$-regular graph, namely $d|U||V|/n$.
 
\begin{lemma}[Expander Mixing Lemma \cite{alon1988explicit}]\label{eml}
For every  $(n,d,\lambda)$-graph $G$ and for every $U, V\subseteq V(G)$ it holds 
 \[\left|e_G(U,V)-\tfrac{d|U||V|}{n}\right|\leq \lambda\sqrt{|U||V|}.\]
\end{lemma}


%% file: 4.1sparse.tex
\section{Immersions in \texorpdfstring{$(n,d,\lambda)$}{(n,d,\textlambda)}-graphs}\label{immersion}
\par In this section, we shall complete the proof of Theorem \ref{ndlthm}. 
According to the value of degree $d$, our constructions of clique immersion are different and we divide the proof of Theorem \ref{ndlthm}  into three cases: sparse case, medium case and dense case as follows.

\begin{lemma}[Sparse case]\label{dhx}
For any $0<\eta<1/2$, there exists $K>0$ such that for sufficiently large $n$ and $d$, every $(n,d,\lambda)$-graph $G$ with $K\lambda<d<\log^{200}n$ contains a $K_{(1-\eta)d}$-immersion.
\end{lemma}

\begin{lemma}[Medium case]\label{dxiao}
There exist positive constants $C$ such that the following holds for all $n$ and $d$ with $\log^{200} n\leq d\leq n/C$. For any $0<\eta<1/2$, every $(n,d,\lambda)$-graph $G$ with $d>2\lambda$ contains a $K_{(1-5\eta)d}$-immersion.
\end{lemma}

\begin{lemma}[Dense case]\label{dda}
For any $c>0$ and $0<\eta<1/2$, there exists $K>0$ such that for sufficiently large $n$, every $(n,d,\lambda)$-graph $G$ with $d\geq cn$ and $d>K\lambda$ contains a $K_{(1-\eta)d}$-immersion.
\end{lemma}

Lemma \ref{dhx} is an immediate corollary from Theorem~\ref{d51} on balanced clique subdivisions. For Lemma \ref{dxiao}, we first verify that $G$ is in fact a robust sublinear expander. Then the uniform distribution of $E(G)$ (see Lemma~\ref{eml}) can be used to find enough edge-disjoint \emph{units}, a notion that was first introduced by Liu, Wang and Yang~\cite{liu2022clique}, which are tree-like structures for building immersions. Finally, we use Lemma \ref{rbsdl} to greedily connect these units, and this process needs to be carried out according to certain rules to achieve that the paths connecting two distinct branch vertices should be edge disjoint.

As to the dense case (Lemma \ref{dda}), we need to find an immersion of a clique of order $d-o(d)$ where $d=\Theta(n)$. In doing this, we adopt another useful method -- R\"odl Nibble. The specific idea is as follows. For any set $A$ of $(1-o(1))d$ vertices in $G$, Lemma \ref{eml} shows that there are roughly $d/n$ portion of edges inside $A$. For the remaining $1-d/n$ portion of  non-adjacent vertex pairs, we first use the R\"odl-Nibble method to connect most of them by pairwise edge-disjoint paths of length two in between $A$ and $V(G)\setminus A$, and finally a small portion of left-over pairs can be greedily connected by edge-disjoint paths of length three (again by Lemma \ref{eml}). By doing the above steps carefully, we can make sure that all the paths are edge disjoint.

\subsection{Embedding immersions in sparse \texorpdfstring{$(n,d,\lambda)$}{(n,d,\textlambda)}-graphs}\label{sparse}
\par Here, we will complete the proof of Lemma \ref{dxiao}. We first give the following lemma.

\begin{lemma}\label{robust}
If $G$ is an $(n,d,\lambda)$-graph with $d\geq 2\lambda$, then $G$ is an $(\epsilon_1,\epsilon_2 d)$-robust-expander for any $0<\epsilon_1\leq 1/8$ and $0<\epsilon_2<1$.
\end{lemma}

\begin{proof}
For every $X\subseteq V(G)$ with $\epsilon_2 d/2\leq |X|\leq n/2$, and any edge set $F\subseteq E(G)$ with $|F|\leq d\cdot\rho(|X|)|X|$, we have 
	\begin{equation}\label{NX}
		|N_{G\setminus F}(X)|\geq \frac{e(X, V(G)\setminus X)-|F|}{d}.
	\end{equation}
By Lemma \ref{Alon}, we have that $e(X, V(G)\setminus X)\geq \frac{(d-\lambda)|X|(n-|X|)}{n}\geq \frac{(d-\lambda)|X|}{2}$. 
And by Inequality (\ref{1m}), we have $\rho(|X|)< \epsilon_1\leq 1/8$ .
So Inequality (\ref{NX}) is at least \[\left(\frac{(d-\lambda)|X|}{2}-d\cdot\rho(|X|)|X|\right)/d\geq \left(\frac{1}{2}-\frac{\lambda}{2d}-\rho(|X|)\right)|X|\geq \rho(|X|)|X|.\]
\end{proof}
The following property shows that after removing some vertices and edges, the average degree of the resulting graph does not decrease too much.
 
\begin{proposition}\label{property}
For $0\leq \eta\leq 1/2$. If $G$ is an $(n,d,\lambda)$-graph, then for any subset $U\subseteq V(G)$ with $|U|\leq \eta n$, and any edge set $W\subseteq E(G)$ with $|W|\leq \frac{\eta(1-\eta)}{2}dn$, we have \[d(G\setminus W-U)\geq (1-3\eta)d.\]
\end{proposition}
	
\begin{proof}
As \[e(G-U)=\frac{nd}{2}-e(G[U])-e(U,V(G)\setminus U)\]
and 
\[e(G[U])+e(U,V(G)\setminus U)\leq |U|d,\]
we have		
\begin{align}\label{G-U}
d(G-U) &\geq \frac{nd-2|U|d}{n-|U|}=d-\frac{|U|d}{n-|U|}\geq d-\frac{\eta}{1-\eta}d\geq d-2\eta d,
\end{align}
where the penultimate inequality holds as $\frac{x}{n-x}$ is an increasing function  by $x$ and the last inequality holds as $0\leq \eta\leq 1/2$. 
Since  
\[|W|\leq \frac{\eta(1-\eta)}{2}dn \leq \frac{1}{2}\eta d (n-|U|),\] we have 
\[d(G\setminus W-U)\geq d(G-U)-\frac{2|W|}{n-|U|}\geq d-2\eta d-\eta d\geq d-3\eta d.\]
\end{proof}
\subsubsection{Units}\label{unit}
\par We first introduce the definition of units from Liu, Wang and Yang \cite{liu2022clique}.
 	
\begin{figure}[H]
	\centering
	\includegraphics[scale=0.16]{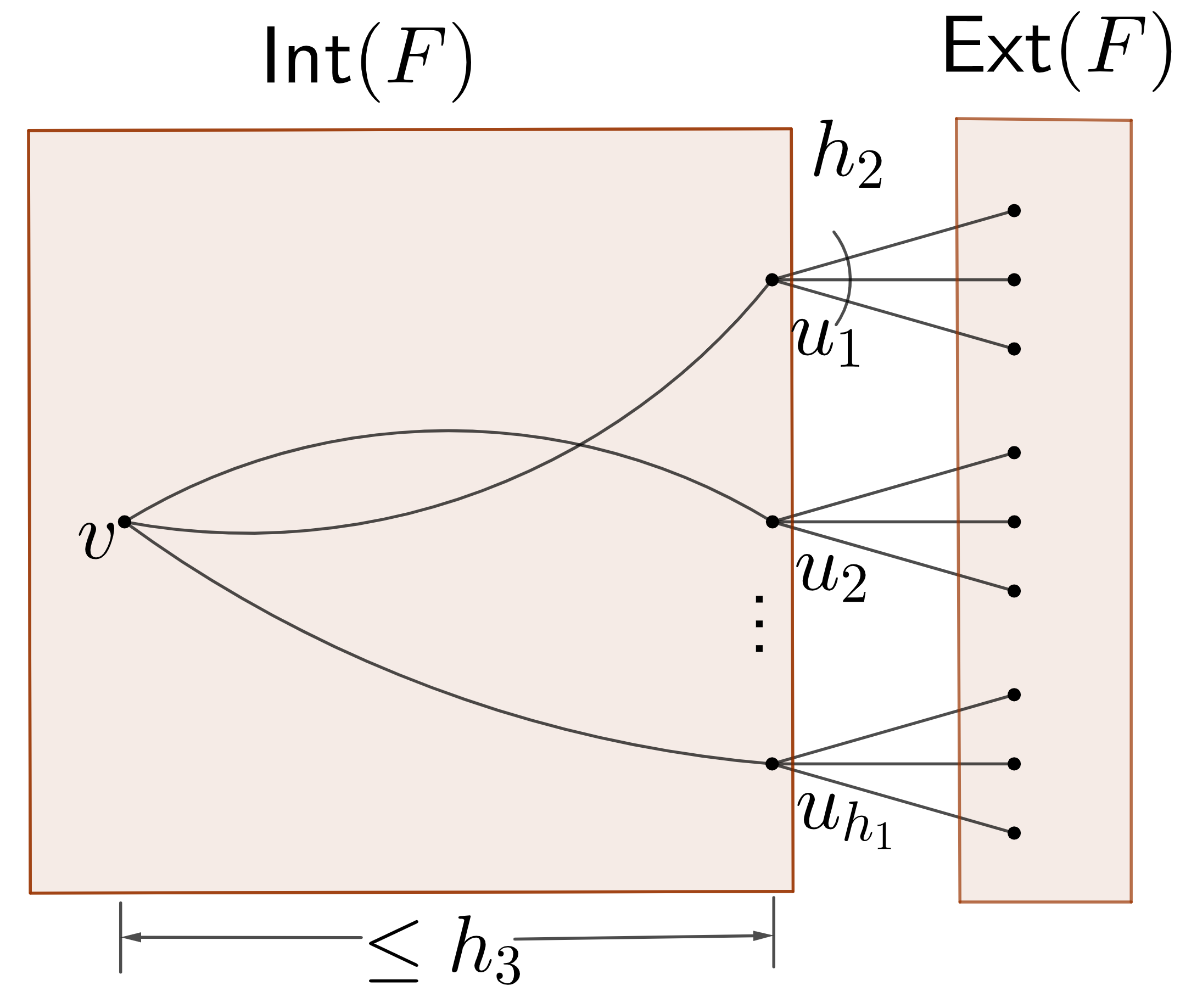}
	\caption{$(h_1, h_2, h_3)$-\emph{unit}: $h_1$ vertex-disjoint stars $S(u_i)$ each of size $h_2$; $h_1$ edge-disjoint $(v, u_i)$-paths each of length at most $h_3$. Note that $u_i$ may appear in a $(v,u_j)$-path for some $i\neq j$.}
	\label{fig1}
\end{figure}

\begin{definition}\rm\cite{liu2022clique}
Given $h_1, h_2, h_3>0$, an \emph{$(h_1, h_2, h_3)$-unit} $F$ is a graph consisting of a center $v$, $h_1$ vertex-disjoint stars $S(u_i)$ centered at $u_i$, each of size $h_2$ and edge-disjoint $(v, u_i)$-paths, for each $i\in[h_1]$, whose length is at most $h_3$. Moreover, the set of interior vertices in all $(v,u_i)$-paths is disjoint from all leaves in $\cup_{i=1}^{h_1} S(u_i)$. The \emph{exterior} of the unit, denoted by $\mathsf{Ext}(F)$, is the set of all leaves in $\cup_{i=1}^{h_1}S(u_i)$, and the rest set of vertices in $F- \mathsf{Ext}(F)$ is called the \emph{interior} of $F$, denoted by $\mathsf{Int}(F)$. We call each $(v,u_i)$-path a \emph{branch} of $F$ and each edge in the star $S(u_i)$ a \emph{pendant edge} (see Figure \ref{fig1}).
\end{definition}

Then we can use Proposition~\ref{property} to greedily find a large number of pairwise edge-disjoint units in an $(n,d,\lambda)$-graph. 
		
\begin{lemma}\label{findunit}
For each $0<\eta<1/2$ and $y\leq 54$, there exists $0<\epsilon_1 \leq 1/8,0<\epsilon_2<1$ and sufficiently large $C>0$ such that the following holds for all $n$ and $d$ with $\log^{200}n\leq d\leq n/C$. If $G$ is an $(n,d,\lambda)$-graph with $d\geq 2\lambda$, then for any vertex set $U\subseteq V(G)$ with $|U|\leq dm^{y+1}/4$, and any edge set $W\subseteq E(G)$ with $|W|\leq d^2 m^{y+1}/3$, $G\setminus W-U$ contains a $(d',m^y,m)$-unit, where $d'=(1-4\eta)d$.
\end{lemma}

The proof of Lemma \ref{findunit} is in the Appendix \ref{appfindunit}. Based on Lemma \ref{findunit}, we can greedily find $d'$ edge-disjoint units in the $(n,d,\lambda)$-graph $G$.
		
\begin{corollary}\label{findunits}
Let $G, d', m, y, n, d, \lambda$ be as in Lemma \ref{findunit}. There are $d'$ edge-disjoint $(d',m^y,m)$-units with distinct centers in $G$. 
\end{corollary}

\begin{proof}
Let $\mathcal{U}$ be a maximum collection of edge-disjoint $(d',m^y,m)$-units in $G$. Suppose for a contradiction that $|\mathcal{U}|<d'$. Let $W,U$ be the set of edges and centers of all units in $\mathcal{U}$, respectively. Then we have
\[|W|<d'(d'm+d'm^y)\leq \frac{d^2m^{y+1}}{3},~\text{and}~~ 
|U|<d.\]
By Lemma \ref{findunit}, $G\setminus W-U$ contains a new $(d',m^y,m)$-unit, a contradiction to the maximality of $\mathcal{U}$.
\end{proof}

\subsubsection{Proof of Lemma \ref{dxiao}}
\par Up to now, we already have that an $(n, d, \lambda)$-graph with $d\geq 2\lambda$ is a robust expander by Lemma \ref{robust}, and 
contains at least $d'$ edge-disjoint units with distinct centers by Corollary \ref{findunits}. The remaining work is to find edge-disjoint short paths connecting any two different centers of $(1-o(1))d$ units.

\begin{proof}[Proof of Lemma \ref{dxiao}]
Let $y:=50$. By Corollary \ref{findunits}, $G$ contains $d'$ edge-disjoint $(d',m^y,m)$-units, say $F_1,\dots,F_{d'}$, whose centers are $v_1, v_2, \dots v_{d'}$, respectively. Let $\mathcal{P}$ be a maximum collection of paths satisfying the following rules. 

\stepcounter{propcounter}
\begin{enumerate}[label=({\bfseries\Alph{propcounter}\arabic{enumi}})]
\item\label{A1} Each path connects $v_i$ and $v_j$ through a $(\mathsf{Ext}(F_i),\mathsf{Ext}(F_j))$-path $P_{ij}$ of length at most $m$, and there is at most one path in $\mathcal{P}$ between any distinct $v_i$ and $v_j$,  $i, j \in [d']$.
\item\label{A2} In each $F_i$, a star is \emph{occupied} if a leaf of it was previously used as an endpoint of a path $P_{ij}$ for some $i\neq j$.
\item\label{A3} Let $v_i, v_j$ be the current pair to connect. Then the $(\mathsf{Ext}(F_i),\mathsf{Ext}(F_j))$-path $P_{ij}$ needs to avoid using
\begin{enumerate}[label=({\alph{propcounter}\arabic*})]
\rm\item \label{b1} any leaf of the occupied stars in $F_i$ or $F_j$ as an endpoint;
\rm\item \label{b2} edges that are in branches of all $F_1,\dots,F_{d'}$ units;
\rm\item \label{b3} edges that are used in previous connecting paths;
\rm\item\label{b4} all centers $v_i$, $i\in[d']$.
\end{enumerate}
\end{enumerate}
By the definition of unit, every path in $\mathcal{P}$ has length at most $4m$.
Then the total number of edges in $\mathcal{P}$ is at most 
\[\binom{d'}{2}\cdot 4m\leq 2d^2 m.\]
The total number of edges in branches of all units is at most 
\[d'\cdot d'\cdot m\leq d^2 m.\] Denote by $Q$ the set of edges we need to avoid in \ref{b2}-\ref{b3}, and then we have $|Q|\leq 3 d^2 m$. In other word, we need to avoid at most $3 d^2 m$ edges and at most $d'$ centers in \ref{A3} in each connection.  
Define a unit to be \emph{bad} if more than $\eta dm^y/2$ pendant edges are used in $\mathcal{P}$, otherwise it is \emph{good}. Then, we get that there are at least $d''=(1-5\eta)d$ good units, say $F_1,\dots,F_{d''}$. Indeed, the number of bad units is at most \[\frac{3d^2m}{\eta dm^y/2}<\eta d/2\leq d'-d''.\]
		
Now we claim that all centers of the units $F_1,\dots,F_{d''}$ are connected via paths in $\mathcal{P}$. Assume for a contradiction that there exists one pair $\{i,j\}\in\binom{[d'']}{2}$ such that there is no path in $\mathcal{P}$ connecting $v_i$ and $v_j$. Let $S_i$ and $S_j$ be the sets of leaves of the occupied stars in $F_i$ and $F_j$, respectively.
Let $F_i'\<F_i\setminus E(\mathcal{P})-S_i+v_i$ and $F_j'\<F_j\setminus E(\mathcal{P})-S_i+v_j$ be two subunits after deletion. Then 
\[|F_i'|\geq (d'-d''-\frac{\eta dm^y/2}{m^y})(m+m^y)\geq \frac{\eta dm^y}{2}=:x,\] and $|F_j'|\geq x$. 
Recall that the number of edges that we need to avoid is at most $3 d^2 m\leq d\cdot \rho(x)\cdot x$. By Lemma \ref{robust}, $G$ is an $(\epsilon_1, \epsilon_2 d)$-robust-expander for $0<\epsilon_1\leq 1/8$ and $1<\epsilon_2<1$. 
Then applying Lemma \ref{rbsdl} with $X_1=F_i', X_2=F_j', W=Q$ and $U=\{v_1,\dots,v_{d'}\}$, there is an $(F_i',F_j')$-path $P_0$ of length at most $m$, satisfying \ref{A3}. Note that $F_i'$ and $F_j'$ are also two units centered at $v_i$ and $v_j$, respectively. So the path $P_0$ can be extended into a $(v_i,v_j)$-path, satisfying \ref{A1} and \ref{A3}. This contradicts the maximality of $\mathcal{P}$.
		
Hence, we can connect all pairs of the units $F_1,\dots,F_{d''}$ with edge-disjoint paths, which implies a desired $K_{d''}$-immersion.
\end{proof}


%% file: 4.2dense.tex


\subsection{Embedding immersions in dense \texorpdfstring{$(n,d,\lambda)$}{(n,d,\textlambda)}-graphs}\label{immersiondense}
\par In this section, we shall finish the proof of Lemma \ref{dda}. 
Since every vertex has degree $d=\Omega(n)$, for any two vertices, one can use Lemma \ref{eml} to show that there are many paths connecting them, but it is not clear how to pick edge-disjoint paths connecting roughly $\binom{d-o(d)}{2}$ pairs. To achieve this, we adopt the approach of R\"odl-Nibble.
In this paper, we use the following version of the R\"odl-Nibble method.

\begin{lemma}{\rm\cite{alon2016probabilistic}}\label{rodl}
For every integer $r\geq 2$ and reals $k\geq 1$ and $\alpha>0$, there are $\gamma=\gamma(r,k,\alpha)>0$ and $d_0=d_0(r,k,a)$ such that for every $n\geq D\geq d_0$ the following holds. Every $r$-uniform hypergraph $H=(V,E)$ on a set $V$ of $n$ vertices in which all vertices have positive degree and which satisfies the following conditions:
\begin{itemize}
\item For all but at most $\gamma n$ vertices $x\in V$  , $d(x)=(1\pm \gamma)D.$ 
\item For all $x\in V$, $d(x)<kD.$
\item For any two distinct $x,y\in V$, $d(x,y)<\gamma D.$	\end{itemize}
contains a matching of at least $(1-\alpha)\frac{n}{r}$ edges.
\end{lemma}

\subsubsection{Regularity of \texorpdfstring{$(n,d,\lambda)$}{(n,d,\textlambda)}-graphs}
\par We first give some notions and tools. Let $G$ be a graph, and let $X,Y\subseteq V(G)$ be disjoint. Then we call $d(X,Y)=e(X,Y)/(|X||Y|)$ the \emph{density} of the pair $(X,Y)$. Given some $\varepsilon>0$, we say a pair $(X,Y)$ of disjoint sets $X,Y\subseteq V$ is \emph{$\varepsilon$-regular} if all $X'\subseteq X$ and $Y'\subseteq Y$ with $|X'|\geq \varepsilon |X|$ and $|Y'|\geq \varepsilon |Y|$ satisfy $|d(X',Y')-d(X,Y)|\leq \varepsilon$.

\begin{lemma}\label{regularity}
For any $c, \delta,\varepsilon>0$, there exists $K\geq 10/(\varepsilon^2 \delta)$ such that the following holds for any sufficiently large $n$. Let $G$ be an $(n,d,\lambda)$-graph with $K\lambda<d=cn$. Then for any disjoint subsets $A,B\subseteq V(G)$ with $|A|\geq |B|\geq \delta d$, $(A,B)$ is an $\varepsilon$-regular pair in $G$ (also $\overline{G}$), and furthermore
\[d_G(A,B)=c\pm\varepsilon, \text{~and~} d_{\overline{G}}(A,B)=(1-c)\pm\varepsilon.\]
\end{lemma}

\begin{proof}
   Since $G$ is an $(n,d,\lambda)$-graph, by Lemma \ref{eml}, we have that for any subsets $A'\subseteq A$ and $B'\subseteq B$ with $|A'|\geq \varepsilon |A|$ and $|B'|\geq \varepsilon |B|$, 
\[d_G(A,B)=\frac{e_G(A,B)}{|A||B|}=\frac{d}{n}\pm \frac{\lambda}{\sqrt{|A||B|}}=c\pm \frac{\lambda}{|B|},\]
\[d_G(A',B')=\frac{e_G(A',B')}{|A'||B'|}=\frac{d}{n}\pm \frac{\lambda}{\sqrt{|A'||B'|}}=c\pm \frac{\lambda}{\varepsilon |B|}.\]
Let $K\geq 10/(\varepsilon^2 \delta)$ and we have
\[|d_G(A,B)-d_G(A',B')|\leq \frac{\lambda}{|B|}+\frac{\lambda}{\varepsilon |B|}<\frac{2\lambda}{\varepsilon |B|}<\frac{2}{\varepsilon \delta K}<\varepsilon.\]
It implies that the pair $(A,B)$ is $\varepsilon$-regular in $G$, and $d_G(A,B)=c\pm\varepsilon.$
Similarly one can verify that the pair $(A,B)$ is $\varepsilon$-regular in $\overline{G}$ with $d_{\overline{G}}(A,B)=(1-c)\pm \varepsilon$.
\end{proof}

We also need the following immediate fact from $\eps$-regularity.

\begin{fact}\label{factc}
Let $d>\varepsilon>0$. If $(I,J)$ is an $\varepsilon$-regular pair of density $d$ and $J'\<J$ has size at least $\eps |J|$, then the number of vertices $u\in I$ satisfying $d(u,J')>(d+\varepsilon)|J'|$ is at most $\varepsilon |I|$, and also the number of vertices $u\in I$ satisfying $d(u, J')<(d-\varepsilon)|J'|$ is at most $\varepsilon |I|$. 
\end{fact}	

\begin{corollary}\label{good}
\rm Let $\varepsilon>0$, $k\in\mathbb{N}$ and pairs $(I, J_i)$ for $i\in [k]$ be $\varepsilon$-regular. Define a vertex $u\in I$ to be \emph{good} for given subsets $J_i'\< J_i$ ($i\in [k]$) with $|J_i'|\ge \eps |J_i|$ if we have $|N(u)\cap J_i|=(d(I,J_i)\pm \varepsilon)|J_i|$, otherwise, it is \emph{bad}. 
Then there are at least $(1-2k\eps)|J|$ good vertices.
\end{corollary}
\subsubsection{Proof of Lemma \ref{dda}}
\par We actually show that for any set of $d-o(d)$ vertices in $G$, we can connect any two distinct vertices of them by edge-disjoint paths. To achieve this goal, we first focus on counting the number of edge-disjoint paths of length two in local structure.

\begin{lemma}\label{oneedge}
For all constants $0<\beta<1$ and $p,q>0$ with $p+q=1$, there exists $\varepsilon=\varepsilon(\beta,p,q)>0$ such that the following holds for every integer $t\geq p/(\varepsilon q)$. Let $k=\lceil qt/p \rceil$. Suppose that $G$ is a tripartite graph with $V(G)=A\cup B\cup C$, where $|A|=|B|=t$ and $|C|=k$. Pairs $(A,B), (A,C)$ and $(B,C)$ are $\varepsilon$-regular with edge densities $q\pm \varepsilon, p\pm\varepsilon,p\pm\varepsilon$ respectively. Then there are at least $(1-\beta)\frac{e(G)}{3}$ edge-disjoint triangles.
\end{lemma}
\begin{proof}
Choose constants $\eps\ll p,q,\beta$.
 We can get the following claim.
\begin{claim}\label{cl1}
Except for at most $\gamma qt^2$ edges, every edge is contained in $(1\pm \gamma)pqt$ triangles, where $\gamma= 12\sqrt{\varepsilon}$.
\end{claim} 

\begin{proof}
By Corollary~\ref{good}, there are at least 
$(1-4\varepsilon)t$ good vertices in $A$ for $B$ and $C$, denoted as $A'\<A$.
Note that for any vertex $a\in A'$,	let $B':=N(a)\cap B$ and $C':=N(a)\cap C$, then we have 
\[|B'|=(q\pm 2\varepsilon)t \text{~and~} |C'|=(p\pm 2\varepsilon)k.\]
Since $(B, C)$ is $\varepsilon$-regular pair and $|B'|\geq \varepsilon |B|$ and $|C'|\geq \varepsilon|C|$, we have $d(B', C')=p\pm 2\varepsilon$.
Similarly, the number of good vertices in $B'$ for $C'$ is at least 
$|B'|-2\varepsilon|B|\ge (1-4\eps/q)|B'|.$
Choose such good vertices for $C'$ from $B'$, denoted as $B''$.
Note that for any vertex $b\in B''$, we have that
$|N(b)\cap C'|=(p\pm 2\varepsilon)|C'|$.
Now we know that the red edge $ab$ is contained in 
$p\pm 2\varepsilon)(p\pm 2\varepsilon)k=(1\pm \gamma)pqt$ triangles, where the last equation follows as $\varepsilon \ll p,q$.
Thus the number of edges between $A$ and $B$ which violate the choices of $ab$ as above is at most 
\begin{align*}
&(|B'|-|B''|)|A'|+(A-|A'|)|B|\\
\leq& 2\varepsilon |B||A'|+4\varepsilon |A||B|\\
\leq&  4\sqrt{\varepsilon} qt^2.
\end{align*}
Thus by the same arguments as above, we can obtain that
in total, we exclude at most $12\sqrt{\varepsilon}qt^2\le \gamma qt^2$ edges.
\end{proof}

Then we construct an auxiliary hypergraph $\mathcal{H}$ as follows.
\stepcounter{propcounter}
\begin{enumerate}[label=({\bfseries\Alph{propcounter}\arabic{enumi}})]	
\item Every edge $e_i$ in $G$ is corresponding to a vertex $v_i$ in $\mathcal{H}$.
\item For three distinct vertices $v_i, v_j, v_k\in V(\mathcal{H})$, $v_i v_j v_k$ is a hyperedge of $\mathcal{H}$ if and only if $e_i e_j e_k$ form a triangle in $G$. 
\end{enumerate}

It is easy to see that the degree of $v_i$ in $\mathcal{H}$ is the number of triangles containing $e_i$ in $G$. The number of edges in $G$ is 
\[e(G)=2(p\pm\varepsilon)tk+(q\pm\varepsilon)t^2\geq 3q(1-\sqrt{\varepsilon})t^2.\]
Then Claim~\ref{cl1} implies that all but at most $\gamma qt^2\leq\gamma e(G)$ vertices $v_i\in V(\mathcal{H})$  have 
\[d_{\mathcal{H}}(v_i)=(1\pm \gamma)pqt.\]
Also each vertex $v_i\in V(\mathcal{H})$ has
 \[d_{\mathcal{H}}(v_i)<\max \{t,k\}<\frac{1}{p^2q} pqt.\]
Moreover, it is trivial that every two distinct $v_i, v_j\in V(\mathcal{H})$ has 
\[d_{\mathcal{H}}(v_i, v_j)\leq 1.\]
Applying Theorem \ref{rodl} with $r=3, k=1/p^2q, D=pqt, n=e(G)$ and the choice of $\eps$, we have that $\mathcal{H}$ contains a matching of at least $(1-\beta)\frac{e(G)}{3}$ edges. It means that $G$ contains at least $(1-\beta)\frac{e(G)}{3}$ edge-disjoint triangles.
\end{proof}

Now, we are ready to finish the proof of Lemma \ref{dda}.

\begin{proof}[Proof of Lemma \ref{dda}]
For convenience, let $d=cn$, $q:=1-c$, $f:=(1-\eta)d$, $t:=\lfloor c\eta^2 d/10\rfloor$, $M_1:=\lfloor f/t \rfloor$, $s:=\lceil qt/c\rceil$, $M_2:=\lfloor(n-f)/s\rfloor$.
We choose $f$ vertices arbitrarily into a set $F$. Then we partition $F$ into $V_0\cup V_1\cup\cdots\cup V_{M_1}$ where $|V_0|<|V_i|=t$ for all $i\in[M_1]$. And the remaining vertices in $G-F$ can be partitioned into $U_0\cup U_1\cup\cdots\cup U_{M_2}$ where $|U_0|<|U_i|=s$ for all $i\in[M_2]$.

We can see that $t\geq c\eta^2d/20$ and $s\geq (1-c)\eta^2d/20$, and let $\delta=\min\{c\eta^2/20,(1-c)\eta^2/20\}$.
Choose $1/K\ll \eps\ll c,\eta$, and we can use Lemma \ref{regularity} to obtain that 
\stepcounter{propcounter}
\begin{enumerate}[label=({\bfseries\Alph{propcounter}\arabic{enumi}})]	
\item\label{C1} any two distinct parts in $V_1,\dots, V_{M_1}, U_1,\dots, U_{M_2}$ are an $\varepsilon$-regular pair in $G$ and $\overline{G}$,
\item for any distinct $i,j\in [M_1]$, $d_G(V_i, V_j)=c\pm\varepsilon$, 
\item for any $i\in [M_1]$ and $j\in[M_2]$, $d_G(V_i, U_j)=c\pm\varepsilon$, 
\item for any distinct $i,j\in [M_1]$, $d_{\overline{G}}(V_i, V_j)=(1-c)\pm\varepsilon$.
\end{enumerate}

If we can prove that any two distinct vertices in $F$ can be connected by pairwise edge-disjoint paths, then we are done. For simplicity, we regard an edge as a path in the following. Firstly, we know that there are some edges in $G[F]$. Avoiding these edges, we want to find edge-disjoint paths to connect the non-adjacent vertices in $F$. We construct an auxiliary graph $G'$ from $G$ as follows (see \autoref{overline(G)}). 
$G'$ has the same vertex set of $G$ and we reserve the partition $\mathcal{Q}:=V_0\cup V_1\cup \cdots \cup V_{M_1}\cup U_0\cup U_1\cup \cdots \cup U_{M_2}$ of $V(G)$ for $V(G')$, and the edges between $F$ and $V(G)\setminus F$ in $G$ and the edges in $E(\overline{G}[F])\setminus \cup_{i=0}^{M_1}E(\overline{G}[V_i])$. We have that in $G'$, any two distinct parts in $V_1,\dots, V_{M_1}, U_1,\dots, U_{M_2}$ are $\varepsilon$-regular, and for any distinct $i,j\in [M_1]$, $d(V_i, V_j)=q\pm\varepsilon$, and for any $i\in [M_1]$ and $j\in[M_2]$, $d(V_i, U_j)=c\pm\varepsilon$. For convenience, we color the edges in $G'[F]$ red, and the other edges black.

We also need the reduced graph $R$ of $G'$ (see \autoref{reduce-bipartite-graph}) where $V(R)=\{v_1,\ldots, v_{M_1}, u_1,\ldots, u_{M_2}\}$ and $v_iv_j$ (or $v_iu_j$) belongs to $E(R)$ if and only if $(V_i, V_j)$ (or $(V_i, U_j)$) is $\varepsilon$-regular. Let $F'=\{v_1,\dots,v_{M_1}\}$. It is easy to see from \ref{C1} that $(F',V(R)\setminus F')$ is a complete bipartite graph and $R[F']$ is a clique of order $M_1$. By Vizing's theorem, we have a decomposition of $E(R[F'])=E_1'\cup\cdots\cup E_{\chi'}'$ where $\chi'\in\{M_1-1, M_1\}$ and each $E_i'$ is a matching of size at least $\frac{M_1-1}{2}$. Thus we can get a partition of $E(G'[F])$ denoted as $\mathcal{W}:=E_0\cup E_1\cup\cdots\cup E_{\chi'}$, where $E_{G'}(V_j,V_k)\subseteq E_i$ if and only if $v_jv_k\in E_i'$ for distinct $j,k\in[M_1]$ and $E_0$ is the set of edges that are incident with $V_0$ or lie inside each $V_i$, $i\in [M_1]$. Note that $|E_0|<tf+M_1\binom{t}{2}$.

\begin{figure} 
\centering
\begin{minipage}{0.45\linewidth}
\parbox[][6cm][c]{\linewidth}{
\centering
\includegraphics[width=0.9\linewidth]{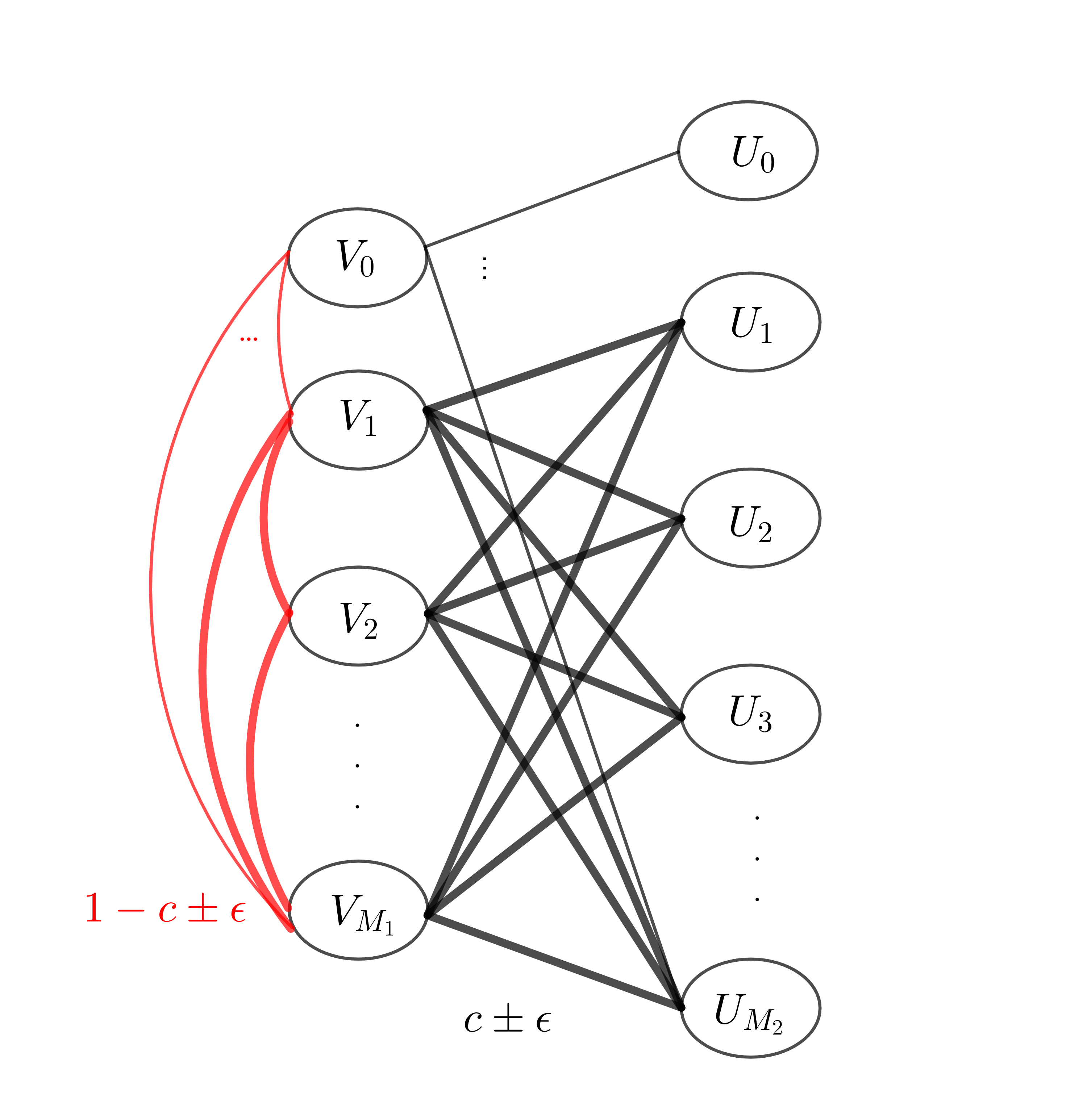}
}
    \caption{The auxiliary graph $G'$}
    \label{overline(G)}
\end{minipage}
\begin{minipage}{0.5\linewidth}
\parbox[][6cm][c]{\linewidth}{
\centering
    \includegraphics[width=0.7\linewidth]{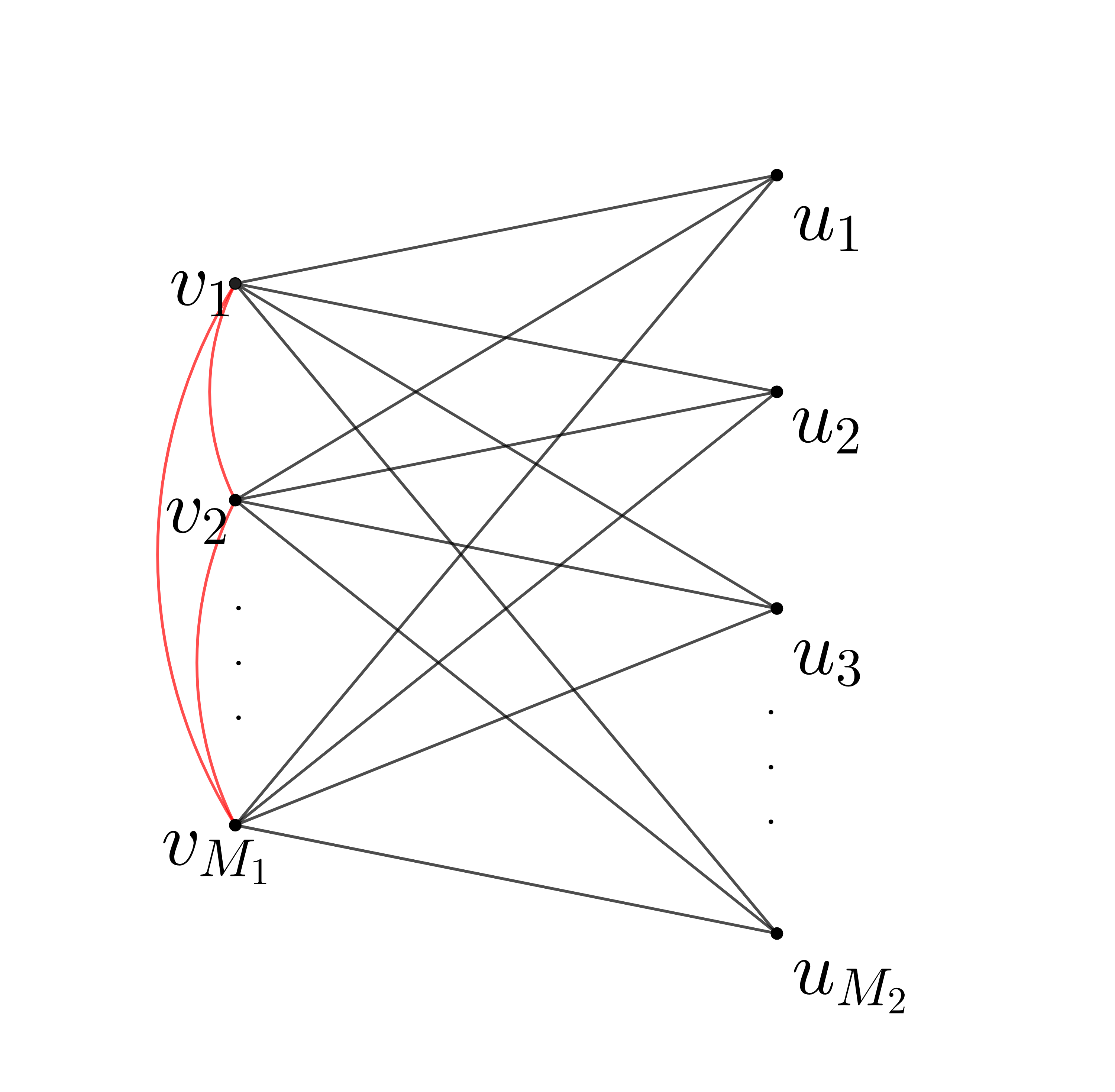}
    }
    \caption{The reduced graph $R$}
    \label{reduce-bipartite-graph}  
\end{minipage}
\end{figure}

We use $2$-\emph{black-path} to denote a $2$-path whose two edges are black.

\begin{claim}\label{2path}
All but at most $c\eta^2 d^2/2$ red edges in $G'[F]$ can be replaced by edge-disjoint $2$-black-paths.
\end{claim}

\begin{proof}
Let $\beta:=c\eta^2/10$. As $n$ is sufficiently large and $M_1\leq f/t$ and \[M_2\geq \frac{n-f}{s}-1\geq \frac{n-f}{qt/c+1}-1\geq \frac{f}{t},\] 
we have $M_2\geq M_1\ge \chi'$. 
For any $j,k\in [M_1]$ and $i\in [M_1]$ with $v_jv_k\in E'_i$, we can see that the induced graph $H:=G'[V_j\cup V_k\cup U_{i}]\setminus (E(V_j)\cup E(V_k)\cup E(U_i))$ is a tripartite graph. Using Lemma \ref{oneedge} on $H$ with $p=c$, $t=t$, $k=s$ and $\beta=\beta$, we can find a suitable $\varepsilon\leq \varepsilon(c,\eta)$ such that there are at least 
\[(1-\beta)\frac{e(H)}{3}=(1-\beta)\frac{(q\pm\varepsilon)t^2+2(c\pm\varepsilon)ts}{3}\geq (1-\beta)(1-\sqrt{\varepsilon})qt^2\geq (1-\beta)^2 qt^2 \] edge-disjoint triangles in $H$.
This implies that all but at most $(q\pm\varepsilon)t^2-(1-\beta)^2 qt^2\leq 3\beta qt^2$ red edges in $E(V_j, V_k)$ can be replaced by edge-disjoint 2-black-paths, denoted as $\mathcal{P}^i_{j,k}$. We repeat the process for every $j,k\in [M_1]$ with $v_jv_k\in E'_i$ and observe that the corresponding families $\mathcal{P}^i_{j,k}$ of paths are pairwise edge disjoint. Therefore we proceed the process on every $i\in[\chi']$ and at the end we obtain a family $\bigcup_{i\in[\chi'],v_jv_k\in E'_i}\mathcal{P}^i_{j,k}$ of edge-disjoint 2-black paths which can replace all but at most 
\begin{equation}\label{uncover1}
    3\beta qt^2\binom{M_1}{2}\leq 2\beta q f^2\leq\frac{1}{5}c\eta^2d^2
\end{equation}
red edges in $\cup_{i=1}^{\chi'} E_i$.

Note that we omit the edges in $E_0$ in the above process, and
\begin{equation}\label{uncover2}
    |E_0|\le tf+M_1\binom{t}{2}\leq \frac{c\eta^2d^2}{10}+\frac{c\eta^2d^2}{20}\leq \frac{3}{20}c\eta^2d^2.
\end{equation}
Combining Inequality (\ref{uncover1}) and (\ref{uncover2}), we know there are at most $c\eta^2 d^2/2$ edges in $\overline{G}[F]$ cannot be replaced by edge-disjoint 2-black-paths.	
\end{proof}
	
Let $\mathcal{P}$ be a maximal collection of paths of length at most two found in the previous two steps. Note that all the edges in $E(\mathcal{P})$ lie inside in $G[F]$ and $G(F,V(G)\setminus F)$. 
So we only need to prove that the remaining at most $c\eta^2 d^2/2$ red edges of $\overline{G}[F]$ can be replaced by edge-disjoint paths avoiding $E(\mathcal{P})$. 
	
If there are two vertices $v_i, v_j \in F$ such that there is no path in $\mathcal{P}$ connecting them, then the number of paths in $\mathcal{P}$ with one end in $v_i$ (or $v_j$) is at most $f-1$. Recall that $G$ is $d$-regular. So there are at least $d-(f-1)\geq \eta d$ edges incident with $v_i$ (or $v_j$) not appearing on the paths in $\mathcal{P}$, and denote the other ends of these incident edges by $N'(v_i)$ and $N'(v_j)$, respectively. Note that $N'(v_i), N'(v_j)\subseteq V(G)\setminus F$. So the edges between $N'(v_i)$ and $N'(v_j)$ are disjoint from $E(\mathcal{P})$. We have
\[|N'(v_i)|\geq \eta d \text{~and~} |N'(v_j)|\geq \eta d.\]
Note that $N'(v_i)\cap N'(v_j)=\emptyset$, otherwise there will be a $(v_i,v_j)$-path of length two. By Lemma \ref{eml} and the choice of $K=K(c,\eta)$, there are at least
	\begin{align*}
	&\frac{d|N'(v_i)||N'(v_j)|}{n}-\lambda\sqrt{|N'(v_i)||N'(v_j)|}\\
		&\geq  c\eta^2 d^2-\eta\lambda d\\
		&\geq \frac{c\eta^2 d^2}{2}
	\end{align*}
 edges between $N'(v_i)$ and $N'(v_j)$. 
This implies that the all the remaining edges $v_iv_j$ of $\overline{G}[F]$ can be greedily linked by edge-disjoint paths of length three. 
Until now, all the vertices of $F$ can be connected to each other. We get a desired $K_{(1-\eta)d}$-immersion.
\end{proof}

%% file: 5balancedimmersion.tex
\section{Balanced immersions in graphs}\label{balancedimmersion}
\par In this section, we will give a proof of Theorem \ref{thm1} which is divided into three cases---dense, sparse, and medium case. The dense case is derived from a result in~\cite{devos2014minimum} on $1$-immersion of cliques, while the sparse case follows a result of Wang on balanced subdivisions in sparse expanders (see Lemma \ref{sparsecase}). The bulk of the work is to solve the medium case, the proof of which will be presented in Subsection \ref{hpgra}, and our main work is to prove the following Lemma \ref{keylemma} and Lemma \ref{dense}, which immediately leads to the medium case of Theorem \ref{thm1}. 

\begin{lemma}\label{keylemma}
There exists $\epsilon_1>0$ such that for every $0<\epsilon_2<1/5$, there are positive constants $K$ such that the following holds for all $n$ and $d$ with $\log^{200}n\leq d\leq n/K$. Let $x,y$ be two positive integers such that $13\leq x<100$ and $y=\lfloor (x-7)/2\rfloor$. If $G$ is an $n$-vertex bipartite  $(dm^x, d^2 m^x, d/2)$-dense $(\epsilon_1, \epsilon_2 d)$-robust-expander with $\delta(G)\geq d$, then $G$ contains a $K_{d/20}^{(\ell)}$-immersion for any odd integer $2\leq \ell\leq m^{y-2}$.
\end{lemma}
\begin{lemma}\label{dense}
 There exist $c>0$ and sufficiently large $K$ such that the following holds for all $n$ and $d$ with $n/d>K$. If $G$ is an $n$-vertex $K_{cd}^{(3)}$-immersion-free graph with $\delta(G)=d$, then $G$ is $(dm^x, d^2 m^x, d/2)$-dense for any integer $0<x<100$.
\end{lemma}
We also use the following result about balanced subdivisions in sparse graphs. 
\begin{lemma}[Sparse case]
{\rm\cite{wang2021balanced}}\label{sparsecase}
There exists $\epsilon_1 > 0$ such that for any $0 < \epsilon_2 < 1/5$ and $s \geq 20$, there exist $d_0 = d_0(\epsilon_1, \epsilon_2, s)$ and some constant $c_1 > 0$ such that the following holds for each $n \geq d \geq d_0$ and $d < \log^s n$. If $G$ is a $K^{(2)}_{d/2}$-subdivision-free bipartite $n$-vertex $(\epsilon_1, \epsilon_2 d)$-expander with $\delta(G) \geq d$, then $G$ contains a $K_{c_1 d}^{(\ell)}$-subdivision for some $\ell\in\mathbb{N}$.
\end{lemma}

Now, we are ready to show Theorem \ref{thm1} based on the above lemmas.

\begin{proof}[Proof of Theorem \ref{thm1}]
There exist constants $0<\epsilon_1,\epsilon_2<1$ such that Corollary \ref{col} holds. Let $d'=d/8$. Then there is a bipartite $(\epsilon_1, \epsilon_2d)$-robust-expander $G'$ in $G$ with $d(G')\geq 2d'$ and $\delta(G')\geq d'$. 
We divide the proof of Theorem \ref{thm1} into three cases
\begin{itemize}
		\item [(i)] When $d'\geq n/K$ for some sufficiently large $K$ as in Lemma \ref{dense}. Applying the result in \cite{devos2014minimum} with $c=1/(16K)$, we can find a $K_{d/(64K)}^{(1)}$-immersion, as desired. 
		\item [(ii)] When $d'\leq \log^{200}n$. Applying Lemma \ref{sparsecase}, we can find a $K_{c' d}^{(\ell)}$-immersion with $c'=\min \{1/16,c_1/8\}$ for some $\ell\in\mathbb{N}$, as desired.
		\item [(iii)] When  $\log^{200}n<d'<n/K$.  Let $c_2$ be a constant playing the role of $c$ in Lemma \ref{dense}. If $G'$ contains a $K_{c_2d'}^{(3)}$-immersion, then we are done. Hence, we may assume that $G'$ is $K_{c_2d}^{(3)}$-immersion-free. Applying  Lemma \ref{dense} with $(x,d)=(10,d')$, we know that $G'$ is $(d'm^{10},d'^2m^{10},d'/2)$-dense. Thus, by Lemma \ref{keylemma}, there is a constant $c_0$ such that $G'$ contains a $K_{c_0 d/8}^{(\ell)}$-immersion for some $\ell\in\mathbb{N}$. 
\end{itemize} 

The proof is completed by taking $c=\min\{1/(64K), 1/16, c_1/8, c_2, c_0/8\}$. 
\end{proof}
\subsection{Reduction to robustly dense graph}
\par
The following lemma is essential to prove Lemma \ref{dense}.

\begin{lemma}\label{drc}
Let $H=(A, B; E)$ be a bipartite graph with edge density $\alpha$, and write $n_1:=|A|$ and $n_2:=|B|$. Then there is a $K_p^{(3)}$-immersion in $H$ for every $p\in\mathbb{N}$ with $p\leq \min\{\alpha n_1/16, \alpha^2 n_2/192\}$.
\end{lemma}
We first give a quick proof of Lemma~\ref{dense}, taking Lemma~\ref{drc} for granted.
\begin{proof}[Proof of Lemma \ref{dense}]
	Let $U\subseteq V(G)$ with $|U|\leq dm^x$ and $W\subseteq E(G)$ with $|W|\leq d^2 m^x$.
	If $|U|<d/4$, then it is trivial that $d(G-U)\geq 3d/4$. Since Inequality (\ref{m200}) and $x<100$, we get $dm^{2x}\leq n/2$. Then  
\[d(G\setminus W-U)\geq d(G-U)-\frac{2|W|}{n-|U|}\geq \frac{3d}{4}-\frac{2d^2m^x}{2dm^{2x}-dm^x}\geq  \frac{3d}{4}-\frac{d}{4}=\frac{d}{2}.\]
	
	Next we assume $d/4\leq|U|\leq dm^x $. Let $G_0:=G\setminus W$. Suppose  for a contradiction that $d(G\setminus W-U)< d/2$. 
	We consider the bipartite graph $G_1:=G_0[U,V(G)-U]$. 
	Since $\delta(G)=d$, $dm^{2x}\leq n/2$, $|U|\leq d m^x$ and $|W|\leq d^2 m^x\leq d|V(G)-U|/8$ as Inequality (\ref{m200}), we have 
	\begin{align*}
		e(G_1)&\geq\sum_{v\in V(G)-U} d_G(v)-2|W|-2e_{G_0}(V(G)-U)\\
		&\geq d|V(G)-U|-\frac{d}{4}|V(G)-U|-\frac{d}{2}|V(G)-U|\geq \frac{d}{4}|V(G)-U|.
	\end{align*}
	
	Recall that the edge density of $G_1$ is $\alpha=e(G_1)/(|U||V(G)-U|)$. Take $c\leq 1/3072$,
	then
	\[\frac{\alpha|U|}{16}=\frac{e(G_1)}{16|V(G)-U|}\geq \frac{d}{64}\geq cd,\]
	and
	\begin{equation}\label{avu}
		\frac{\alpha^2 |V(G)-U|}{192}\geq \frac{d^2|V(G)-U|}{3072|U|^2}>\frac{|V(G)-U|}{3072m^{2x}}.
	\end{equation}
As $0<x<100$ and $dm^{2x}\leq n/2$, we get $|V(G)-U|>n/2$. Thus, the Inequality (\ref{avu}) implies that 
	\[\frac{\alpha^2 |V(G)-U|}{192}>\frac{d}{3072}\geq cd.\]
	
	Hence, applying Lemma \ref{drc} with $(U,V(G)-U;E(G_1))=(A,B;E)$, we can find a $K_{cd}^{(3)}$-immersion in $G$, a contradiction.
\end{proof}

\begin{proof}[Proof of Lemma \ref{drc}]
		Let $p$ be an integer satisfying $p\leq\min\{\alpha n_1/16, \alpha^2 n_2/192\}$.
		We pick a vertex $b$ from $B$ uniformly at random. Let $A_0:=N_H (b)$ and the random variable $X:=|A_0|$. Then
		\[\mathbb{E}[X]=\sum_{v\in A} \frac{d(v)}{n_2}=\frac{e(H)}{n_2}=\frac{\alpha n_1 n_2}{n_2}=\alpha n_1.\]

		We say a pair $(u,v)\in A$ is \emph{bad} if $d(u,v)<3p$ .
		Let $Y$ be the random variable counting the number of bad pairs in $A_0$. Then
		\[\mathbb{E}[Y]\leq \binom{n_1}{2}\cdot \frac{3p}{n_2}\leq \frac{3 p n_1 ^2}{2 n_2}.\]
		
		Since $\mathbb{D}[X]=\mathbb{E}[X^2]-\mathbb{E}[X]^2\geq 0$, we have $\mathbb{E}[X^2]\geq \mathbb{E}[X]^2$. Then by the linearity of expectation, we obtain  
		\[\mathbb{E}\left[X^2-\frac{\mathbb{E}[X]^2 Y}{2\mathbb{E}[Y]}-\frac{\mathbb{E}[X]^2}{2} \right]\geq 0.\]
		So there is a choice of $b\in B$ such that 
		\begin{equation*}
			X^2\geq \frac{\mathbb{E}[X]^2}{2}>\frac{\alpha^2 n_1^2}{2} \ \ \text{and} \ \  Y\leq \frac{2\mathbb{E}[Y]X^2}{\mathbb{E}[X]^2}\leq \frac{3pX^2}{\alpha^2 n_2}.
		\end{equation*}
		It implies that we can find a set $A_0$ with $|A_0|=X >\alpha n_1/\sqrt{2}>\alpha n_1/2$ and the number of bad pairs in $A_0$ is at most $(3p|A_0|^2)/(\alpha^2 n_2)$.

		Let $A_1$ be a subset of $A_0$ in which every vertex has codegree less than $3p$ with more than $|A_0|/16$ other vertices of $A_0$. Since $p\leq \alpha^2 n_2/192$, we have
		\[|A_1|\leq \frac{32Y}{|A_0|}\leq\frac{96p|A_0|}{\alpha^2 n_2}\leq\frac{|A_0|}{2}.\]
		
		Let $A_2:=A_0-A_1$. We have $|A_2|\geq |A_0|/2>\alpha n_1/4\geq 4p$ since $p\leq \alpha n_1/16$. 
		Let $A_3, A_4$ be a bipartition of $A_2$ with $|A_4|=p$.
		
		We claim that for each pair of vertices in $A_4$, we can find a path of length $4$ to connect them such that all  these paths are pairwise edge disjoint, which yields a desired $K_p^{(3)}$-immersion. Assume that $u,v$ is the current pair in  $A_4$  to be connected. Recall that the number of bad pairs in $A_2$ containing $u$ or $v$ is at most $|A_0|/8$. So there are at least $|A_2|-|A_4|-|A_0|/8 \geq |A_0|/8>p$ vertices $a_i$, say $A_3'$, in $A_3$ such that the codegree $d(u,a_i)$ and $d(v,a_i)$ are both at least $3p$. 
  For the vertex $u$, there are at most $2(p-1)$ vertices in $N_H(u)$ appearing on the previous connections with one endpoint $u$, and similarly it holds for $v$. We say $a_i\in A_3'$ is a \emph{nice} vertex if the number of vertices in $N_H(a_i)$ appearing on all previous connections is at most $p$. Note that connecting every pair of vertices in $A_4$ only needs $\binom{p}{2}$ edge-disjoint paths, and one such path occupies two vertices in the neighbor of some vertex in $A_3'$. So the number of nice vertices in $ A_3'$ is at least 
		\[|A'_3|-\frac{2\binom{p}{2}}{p}\geq p-(p-1)=1.\]
        This implies that there exists a nice vertex in $A_3'$, say $a$. 
		 Finally, for $u$, there are at least $d(u,a)-2(p-1)-p\geq 3p-2(p-1)-p\geq 2$ vertices in $N_H(u,a)$ which do not appear on the previous connections with one endpoint $u$, and similarly for $v$. So we can choose two distinct vertices $b_i\in N_H(u,a), b_j\in N_H(v,a)$ to form a desired path $ub_iab_jv$.
		
		As above, we can greedily connect every pair of vertices in $A_4$ via pairwise edge-disjoint paths in $H$.
	\end{proof}

\subsection{The medium case}\label{hpgra}
\par In this section, we shall complete the proof of Lemma \ref{keylemma}. First, we need to find many edge-disjoint units. In order to form a balanced immersion, we shall construct a useful gadget called \emph{adjuster}, introduced by Liu and Montgomery\cite{liu2020solution}  to control the length of the paths. Second, we use robust small diameter lemma to connect centers of distinct units greedily with adjusters. By doing this step carefully, we ensure that the paths connecting different centers are edge disjoint and have the same length. 
In the following subsections, we will introduce some useful concepts from \cite{liu2020solution,liu2022clique} and some powerful tools. 

\subsubsection{Finding units}
\par
The following lemma guarantees the existence of a desired unit after avoiding some vertices and edges, and we postpone its proof in the Appendix \ref{appfu}.

\begin{lemma}\label{fu}
	For each $0<\epsilon_1,\epsilon_2<1$, there exists $K>0$ such that the following holds for all $n$ and $d$ with $\log^{200} n\leq d\leq n/K$. Let $x,y$ be two positive integers such that  $x<100$ and $y=\lfloor (x-7)/2\rfloor$. If $G$ is an $n$-vertex $(dm^x, d^2 m^x, d/2)$-dense $(\epsilon_1,\epsilon_2 d)$-robust-expander with $\delta(G)=d$, 
	then for any subset $U\subseteq V(G)$ with $|U|\leq dm^{y+1}$, and any $W\subseteq E(G)$ with $|W|\leq d^2 m^y$, we can find a $(3d/20,m^y,m)$-unit in $G\setminus W-U$.
\end{lemma}

We immediately get a large collection of desired units with distinct centers.

\begin{corollary}\label{findunitcor}
	Let $\epsilon_1,\epsilon_2, K, n, d, x, y$ be as in Lemma \ref{fu}. If $G$ is an $n$-vertex $(dm^x, d^2 m^x, d/2)$-dense $(\epsilon_1,\epsilon_2 d)$-robust-expander with $\delta(G)=d$, then there are at least $d/5$ edge-disjoint $(3d/20,m^y,m)$-units, say $F_1,\dots,F_{d/5}$, with distinct centers $v_1,\dots,v_{d/5}$. 	
\end{corollary}
\begin{proof}
	Suppose for a contradiction that $G$ contains less than $d/5$ edge-disjoint $(3d/20,m^y,m)$-units. Let $W$ and $U$ be the set of edges and the set of centers of these $(3d/20,m^y,m)$-units, respectively. Then we have
	$$|W|<\frac{d}{5}\cdot\left(\frac{3d}{20}\cdot m+\frac{3d}{20}\cdot m^y\right)\leq d^2m^y,$$
	and 
	$$|U|<\frac{d}{5}\leq dm^{y+1}.$$
	By Lemma \ref{fu}, there exists one more $(3d/20,m^y,m)$-unit in $G\setminus W-U$, a contradiction.
\end{proof}
\subsubsection{Building adjusters}

\begin{definition}{\rm\cite{liu2020solution}}
\rm	 Given a vertex $v$ in a graph $F$, $F$ is a $(D,m)$-\emph{expansion} of $v$ if $|F|=D$ and for any vertex $u\in V(F)\setminus \{v\}$, the distance between $u$ and $v$, denoted by $dist(u,v)$, is at most $m$.
\end{definition}

By the definition of $(D,m)$-expansion, it is easy to get the following property.

\begin{proposition}{\rm\cite{liu2020solution}}\label{DD}
	Let $D,m\in \mathbb{N}$ and $1\leq D'\leq D$. Then any graph $F$ which is a $(D,m)$-expansion of $v$ contains a subgraph which is a $(D',m)$-expansion of $v$.
\end{proposition}

Now we give the definition of \emph{adjuster}, which was first introduced by Liu and Montgomery \cite{liu2020solution}.

\begin{definition}
	\rm A \emph{$(D,m,k)$-adjuster} $\mathcal{A}=(u_1,I_1,u_2,I_2,A)$ in a graph $G$ consists of \emph{core vertices} $u_1,u_2\in V(G)$, graphs $I_1,I_2\subseteq G$ and a \emph{center vertex set} $A$ such that the following holds for some $\ell\in \mathbb{N}$.
	\stepcounter{propcounter}
	\begin{enumerate}[label = ({\bfseries \Alph{propcounter}\arabic{enumi}})]
		\rm\item\label{adjla1} $A$, $V(I_1)$ and $V(I_2)$ are pairwise disjoint.
		\rm\item\label{adjla2} For each $i\in[2]$, $I_i$ is a $(D,m)$-expansion of $u_i$.
		\rm\item\label{adjla3} $|A|\leq 10mk$.
		\rm\item\label{adjla4} For each $i\in\{0,1,\ldots,k\}$, there is a $(u_1,u_2)$-path in $G[A\cup\{u_1,u_2\}]$ with length $\ell+2i$.
	\end{enumerate}
\end{definition}

The graphs $I_1$ and $I_2$ are the \emph{ends} of the adjuster, and  we have $V(\mathcal{A})=V(I_1)\cup V(I_2)\cup A$. We denote by $\ell(A)$ the smallest $\ell$ such that \ref{adjla4} holds. Note that $\ell(A) \leq |A| + 1 \leq 10mk + 1.$
The following lemma guarantees the existence of an adjuster after deleting a considerable number of vertices and edges simultaneously. 

\begin{lemma}\label{adjuster}
	For every $0<\epsilon_1,\epsilon_2<1$, there exists $K>0$ such that the following holds for all $n$ and $d$ with $\log^{200} n\leq d\leq n/K$. Let $x, y$ be two positive integers such that $y<x/2-3<28$ and $D=dm^y$. 
	If $G$ is an $n$-vertex $(dm^x, d^2 m^x, d/2)$-dense $(\epsilon_1,\epsilon_2 d)$-robust-expander graph with $\delta(G)=d$. Then for any subset $U\subseteq V(G)$ with $|U|\leq D/(10m)$, and any $W\subseteq E(G)$ with $|W|\leq d^2 m^{y-1}$, $G\setminus W-U$ contains a $(D, m, r)$-adjuster for any $r\leq dm^{y-3}$.
\end{lemma}

The idea to prove Lemma \ref{adjuster} is similar as Lemma 3.9 in \cite{Luan2022BalancedSO} and so we present its proof in Appendix \ref{appadjuster}.
\subsubsection{Connecting vertices with edge-disjoint paths}

\begin{lemma}\label{fp}
	For every $0<\epsilon_1,\epsilon_2<1$, there exists $K>0$ such that the following holds for all $n$ and $d$ with $\log^{200}n\leq d\leq n/K$. Let integers $x<200$, $y=\lfloor (x-7)/2\rfloor$ and $D\geq 4dm^{y-1}$. 
	Suppose that $G$ is an $n$-vertex $(dm^x, d^2 m^x, d/2)$-dense $(\epsilon_1,\epsilon_2 d)$-robust-expander graph with $\delta(G)=d$. 
	Let $Z_1, Z_2\subseteq V(G)$ be vertex sets of size at least $D$, and $I_1, I_2\subseteq V(G)$ be vertex-disjoint $(D,m)$-expansions of $v_1$ and $v_2$, respectively.
	Then for any $\ell\leq dm^{y-2}$, any subset $U\subseteq V(G)$ with $|U|\leq dm^{y-2}$, and any $W\subseteq E(G)$ with $|W|\leq d^2 m^{y-2}$, $G\setminus W-U$ contains vertex-disjoint paths $P$ and $Q$ with $\ell\leq\ell(P)+\ell(Q)\leq\ell+10m$ such that $P$,$Q$ link $\{z_1,z_2\}$ to $\{v_1,v_2\}$ for $z_i\in Z_i$ and $z_1\neq z_2$.
\end{lemma}

We postpone the proof of Lemma \ref{fp} in Appendix \ref{appfp}. And we are ready to prove Lemma \ref{keylemma}.

\begin{proof}[Proof of Lemma \ref{keylemma}]
By Corollary \ref{findunitcor}, there are at least $d/5$ edge-disjoint $(3d/20,m^y,m)$ units $F_1,\dots,F_{d/5}$ centered at $v_1, \dots, v_{d/5}$ in $G$. 	
Thus, there are at least $d/10$ units among  $F_1,\dots,F_{d/5}$ such that their centers are in the same part of the bipartition of $G$, say  $F_1,\dots,F_{d/10}$.

Now we aim to construct a balanced immersion by  $F_1,\dots,F_{d/10}$.
Let $\mathcal{P}$ be the maximum collection of paths under the following rules.
	\stepcounter{propcounter}
\begin{enumerate}[label = ({\bfseries \Alph{propcounter}\arabic{enumi}})]
	\rm\item\label{rule1} Each path connects $v_i$ and $v_j$ via an $(\mathsf{Ext}(F_i),\mathsf{Ext}(F_j))$-path of length at most $m$, and there is only one path in $\mathcal{P}$ between $v_i$ and $v_j$, for $i, j \in [d/10]$ with $i\neq j$.
	\rm\item\label{rule2} All paths in $\mathcal{P}$ are pairwise edge-disjoint and have the same length $\ell+1$.
   \rm\item\label{rule4} In each $F_i$, a star is \emph{occupied} if a leaf of it was previously used as an endpoint in an $(\mathsf{Ext}(F_i),\mathsf{Ext}(F_k))$-path for some $i\neq k$.
	\rm\item\label{rule3} Let $v_i, v_j$ be the current pair to connect. Then the subpath between $\mathsf{Ext}(F_i)$ and $\mathsf{Ext}(F_j)$ of $(v_i,v_j)$-path needs to avoid 
 \begin{itemize}
     \item any leaf of the occupied stars in $F_i$ or $F_j$ as an endpoint.
     \item edges that are in branches of $F_1,\dots,F_{d/10}$;
     \item edges that are used in previous connections;
     \item vertices that are in branches of $F_i$ or $F_j$;
     \item all centers $v_1,\dots,v_{d/10}$;
 \end{itemize}
\end{enumerate}
To be convenient, we denote by $E(\mathcal{P})$ the set of edges in $\mathcal{P}$.
Since $\ell\leq m^{y-2}$, we have  $|E(\mathcal{P})|\leq\binom{d/10}{2}(\ell+1)\leq d^2 m^{y-2}/2$. Let $W_1$ be the set of all edges of the branches in $F_1,\dots,F_{d/10}$. Then   $|W_1|\leq 3d/20\cdot m\cdot d/10\leq d^2 m$.
A unit $F$ is \emph{bad} if more than $3d/20\cdot m^{y-1}$ pendant edges of $F$ are used in $\mathcal{P}$. Note that the number of bad units is at most $\frac{d^2m^{y-2}/2}{3d/20\cdot m^{y-1}}<d/20$. So there are at least $d/10-d/20=d/20$ units which are not bad, say $F_1,\dots,F_{d/20}$.

The following claim implies Lemma  \ref{keylemma} holds.
\begin{claim}\label{111}
	There are pairwise edge-disjoint paths of length $\ell+1$ for every pair $v_i, v_j$, $i, j \in [d/20]$ with $i\neq j$.
\end{claim}
\begin{proof}[Proof of Claim \ref{111}]
	For a contradiction, we assume that there exist $v_i$ and $v_j$ such that no path in $\mathcal{P}$ connects them.
	Denote by $U_1$ the set of vertices in branches of $F_i$ or $F_j$, then $|U_1|\leq 2\cdot 3d/20\cdot m\leq dm$.
	Let $W_2=W_1\cup E(\mathcal{P})$, then $|W_2|\leq|W_1|+|E(\mathcal{P})|\leq d^2 m^{y-2}.$ 
	Applying Lemma \ref{adjuster} with $U=U_1$ and $W=W_2$, there is a $(dm^y,m,20m)$-adjuster $\mathcal{A}=(u_1,I_1,u_2,I_2,A)$ in $G\setminus W_2-U_1$. Note that $|A|\leq 200m^2$ and $\ell(A)\leq |A|+1\leq 210m^2$.

 Let $F'_i$ and $F'_j$ be the remaining connected subgraph of $F_i$ and $F_j$ after deleting the edges in $\mathcal{P}$. Note that the number of paths in $\mathcal{P}$ with one end in $\{v_i,v_j\}$ is at most $d/20$. Since $F_i$ and $F_j$ are not bad units, for each $k\in \{i,j\}$, $F_k$ has at most $3 d/20\cdot m^{y-1}/m^y<3d/(20m)$ stars whose edges are all used in $\mathcal{P}$. Recall that $F_i$ and $F_j$ are $(3d/20, m^y, m)$-units, so for each $k\in \{i,j\}$, we can preserve at least $d/20$ branches and stars of $F_k$ for $F'_k$. Then we have $|F'_i|$ and $|F'_j|$ are both at least $d/20\cdot m+d/20\cdot m^y-3d/20\cdot m^{y-1}\geq 4dm^{y-1}$. Let $\ell'=\ell+1-14m-\ell(\mathcal{A})$.
	Applying Lemma \ref{fp} with $(Z_1,Z_2,I_1,I_2,U,W)=(F'_i,F'_j,I_1,I_2,A\cup U_1,W_2$, we can find vertex-disjoint $(F_i,u_1)$-path $P_1$ and $(F_j,u_2)$-path $Q_1$ with $\ell'\leq \ell(P_1)+\ell(Q_1)\leq \ell'+10m$. It is easy to know that $P_1$ and $Q_1$ can be extended into vertex-disjoint $(v_i,u_1)$-path $P$ and $(v_j,u_2)$-path $Q$, respectively, such that $\ell'\leq \ell(P)+\ell(Q)\leq \ell'+10m+2m+2\leq \ell'+14m$. Then we have $\ell(\mathcal{A})\leq\ell+1-\ell(P)-\ell(Q)\leq \ell(\mathcal{A})+14m$. By the property of the adjuster,  there is a $(u_1,u_2)$-path $R$ of length $\ell+1-\ell(P)-\ell(Q)$ inside $A$. Thus we find a $(v_i,v_j)$-path $v_i P u_1 R u_2 Q v_j$ satisfying \ref{rule1}-\ref{rule3}, a contradiction.
\end{proof}
In conclusion, we complete the proof of Lemma \ref{keylemma}.
\end{proof}

%% file: 6subdivision.tex
\section{Balanced subdivisions in \texorpdfstring{$(n,d,\lambda)$}{(n,d,\textlambda)}-graphs}\label{subdivision}
\par In this section, we aim to finish  the proof of Theorem \ref{d21}. By Lemma \ref{d51}, Dragani{\'c}, Krivelevich and Nenadov~\cite{draganic2022rolling} solved the case when $d\leq n^{1/5}/2$. Based on their result, we prove the case when $n^{1/5}/2<d <n^{1/2-o(1)}$ and the proof combines a random sample argument with an application of the rolling-back method in \cite{draganic2022rolling} for embedding disjoint paths of a fixed length. First, we give several useful results. The first one is a simple lemma for greedily embedding stars in $d$-regular graphs.

\begin{lemma}\label{star}
For any constant $0\leq \eta\leq 1/4$, every $d$-regular graph $G$ with $2/\eta\leq d\leq\eta\sqrt{n}$ has $d$ disjoint copies of star $K_{1,(1-\eta/2)d}$.
\end{lemma}

\begin{proof}
Suppose that there are only $t<d$ disjoint copies of $K_{1,(1-\eta/2)d}$ in $G$. Let $U$ be the set of vertices of all such stars. Then 
\[|U|=t+(1-\frac{\eta}{2})dt\leq\eta\sqrt{n}+(1-\frac{\eta}{2})\eta^2 n\leq\eta^2 n.\]
By Inequality (\ref{G-U}), 
\[d(G-U)\geq (1-2\eta^2)d\geq(1-\frac{\eta}{2})d.\]
This implies that we can find one more $K_{1,(1-\eta/2)d}$ in $G-U$, which is a contradiction.
\end{proof}

 The following definition and lemmas in~\cite{Dellamonica2008An,draganic2022rolling} are powerful tools to connect leaves of different stars with paths of equal length. 
\begin{definition}{\rm\cite{draganic2022rolling}}
Let \( G = (V, E) \) be a graph and $\partial_G(x)$ be the set of edges incident with vertex $x$ in $G$. We say that $G$ has property \( P_{\alpha}(n, d) \) if for every \( X \subseteq V \) of size \( |X| \leq n \) and every \( F \subseteq E \) such that \( |F \cap \partial_G(x)| \leq \alpha \cdot d_G(x) \) for every \( x \in X \), we have \( |N_{G - F}(X)| \geq 2d|X| \).

\end{definition}

\begin{lemma}{\rm\cite{Dellamonica2008An}}\label{proper}
	Let $G$ be an $(n, d, \lambda)$-graph and $d_0$, $n_0$ be positive integers.
	$G$ has \emph{property $P_{\alpha}(n_0, d_0)$} for $\alpha>0$ if the following holds:
	\[1-\alpha>\frac{n_0(1+4d_0)}{2n}+\frac{\lambda}{d}(1+\sqrt{2d_0}).\]	
\end{lemma}

\begin{lemma}{\rm\cite{draganic2022rolling}}\label{connect}
	Let $G$ be a graph with the $P_{\alpha}(n_0,d_0)$ property for $3\leq d_0<n_0$, and such that for every two disjoint $U,V\subset V(G)$ of size $|U|,|V|\geq n_0(d_0-1)/16$, there exists an edge between $U$ and $V$. Let $S'$ be any set of vertices such that $|N_G(x)\cap S'|\leq\beta d_G(x)$ for every $x\in V(G)$ and let $P=\{a_i,b_i\}$ be a collection of at most $\frac{d_0n_0\log{d_0}}{15\log{n_0}}$ disjoint pairs from $S'$. If $\beta\leq 2\alpha-1$, then there exist vertex-disjoint paths in $G$ between every pair of vertices ${a_i,b_i}$ such that the length of each path is $2 \left\lceil \frac{\log(n_0/16)}{\log(d_0-1)} \right\rceil + 3$.
\end{lemma}

\begin{proof}[Proof of Theorem \ref{d21}]
Suppose $d\leq n^{1/5}/2$. Simply let $d_0=3$. We have $240\lambda < d \leq n^{1/5}/2$ and $t=d-80\lambda\sqrt{d_0}>(1-\eta)d$, which finishes the proof by using Theorem \ref{d51}.
	
 Suppose $n^{1/5}/2<d<n^{1/2-\epsilon}$. By Lemma \ref{star}, there are $t:=(1-\eta)d$ disjoint copies of $K_{1,(1-\eta/2)d}$ in $G$, say $S_1,\dots,S_{t}$ with centers $u_1,\dots, u_t$ respectively. Let $U=\{u_1,\ldots,u_t\}$ and $A_i$ be the set of all leaves in $S_i$.
Now we construct a random set $S\subset V(G)\backslash U$, where each vertex in $V(G)\backslash U$ is independently chosen to $S$ with probability $p=1-\eta/4$. Since $A_i\cap U=\emptyset$, we have $|A_i\cap S|\sim Bin(|A_i|,p)$. Clearly $\mu:=\mathbb{E}[|A_i\cap S|]=(1-\eta/2)(1-\eta/4)d$. Using Chernoff inequality, we have
	\[\mathbb{P}[|A_i\cap S|<(1-\eta)d]<\exp(-\tfrac{\mu\eta^2(2+\eta)^2}{2(2-\eta)^2(4-\eta)^2})<\exp(-\tfrac{\eta^2d}{32})<\exp(-\tfrac{\eta^2 n^{1/5}}{64}).\]
	Therefore,
 \[\mathbb{P}[\bigcup_{i\in[d]} \{|A_i\cap S|<(1- \eta )d\}]<n\exp(-\tfrac{\eta^2 n^{1/5}}{64}).\]
	For every vertex $v\in V(G)$, we have $|N(v)\backslash (U\cup S)|\sim Bin(|N(v)\backslash U|,1-p)$. 
 Since $|U|=(1-\eta)d$, $|N(v)\backslash U|\geq \eta d$, we have $\mathbb{E}[|N(v)\backslash(U\cup S)|]=|N(v)\backslash U|(1-p)=\eta^2 d/4$.
 Using Chernoff inequality, we have
 \[\mathbb{P}[|N(v)\backslash (U\cup S)|<\eta^2d/8]<\exp(-(\eta^2d/4)/8)<\exp(-\tfrac{\eta^2n^{1/5}}{64}).\]
	By taking a union bound, we have
	\[\mathbb{P}[\bigcup_{v\in V(G)}\{ |N(v)\backslash (U\cup S)|<\eta^2d/8\}]<n\exp(-\tfrac{\eta^2 n^{1/5}}{64}).\]
	Notice that there is $N_0(\eta)>0$ such that $n\exp(-\frac{n^{1/5}\eta^2}{64})<\frac{1}{2}$ holds for $n>N_0(\eta)$. 
 Therefore, with positive probability, $|N(v)\backslash (U\cup S)|\geq\eta^2d/8$ holds for every $v\in V(G)$ and $|A_i\cap S|\geq(1-\eta )d$ holds for every $i\in [t]$. Denote the random set satisfying all these conditions by $S_0$.
	
	Now we can embed a subdivision of $K_t$ with the help of Lemma \ref{connect}. Let $u_i$ be the $i$th branching vertex of the subdivision. Since $|A_i\cap S_0|\geq t$, we can arbitrarily choose different leaves $u_{i1}$,\dots,$u_{it}$ from $A_i\cap S_0$ for every $i\in[t]$. Each of these vertices is connected to $u_i$ in $G$. Let $S'=U\cup(\cup_{i,j\in[t]}\{u_{ij}\})$. Since $|N(v)\backslash(U\cup S_0)|\geq\eta^2d/8$ holds for every $v\in V(G)$, $|N(v)\cap S'|\leq (1-\eta^2/8)d$. Set \[\alpha=1-\eta^2/16,~~ \beta=2\alpha -1=1-\eta^2/8,~~ d_0=3,~  n_0=\eta^2n/256.\] As $d\geq 240\lambda/\eta^2$,  
by Lemma \ref{proper}, $G$ has property $P_{\alpha}(n_0,d_0)$. Notice that $ |N(v) \cap S'|\leq \beta d(v)$ and 
	\[\frac{d_0n_0\log(d_0)}{15\log(n_0)}>\frac{\eta^2n}{2048\log(n)}>(1-\eta)^2 d^2\]
	holds for $n>N_1(\eta)$ with some sufficiently large $N_1(\eta)$. Also, for every disjoint $U,V\subset V(G)$, $|U|,|V|\geq n_0(d_0-1)/16=\eta ^2n/2048$, 
 by Lemma \ref{eml}, 
	\[e(U,V)\geq \frac{d|U||V|}{n}-\lambda\sqrt{|V||U|}>0,\]
	which means that there is at least one edge between $U$ and $V$. Therefore, $G$ and  $S'$ satisfy all conditions in Lemma \ref{connect}. Using Lemma \ref{connect}, we can connect every pairs of $u_{ij}$ and $u_{ji}$ for $i,j\in[t], i\neq j$ with disjoint paths. Together with the edges between $U$ and $\cup_{i,j\in [t]} \{u_{ij}\}$, we embeds a subdivision of $K_t$ into $G$, where the length of each path is $2 \left\lceil \frac{\log(n_0/16)}{\log(d_0-1)} \right\rceil + 5=2 \left\lceil \log(\eta^2n/4096) \right\rceil + 5$, as desired.
\end{proof}

\noindent
\textbf{Remark.} Notice that by adding extra constraint on $\lambda$, we may work out the following version of Theorem \ref{d21}, which matches the best possible bound $d=O(n^{1/2})$.
\begin{theorem}
    For every $0<\eta<1/2$, there is $\epsilon>0$ such that every $(n,d,\lambda)$-graph $G$ with $\frac{1}{\epsilon}<d<\epsilon n^{1/2}$, $\lambda<d^{1-\eta}$ contains a subdivision of $K_t$ where $t=(1-\eta)d$.
\end{theorem}
The proof is simply repeating the above proof of Theorem \ref{d21}, but setting $d_0=n^{\eta}$ and $n_0=\frac{\eta}{8}n^{1-\eta}$. It is easy to verify that $G$ has $P_{\alpha}(n_0,d_0)$ property and the statement follows by using Lemma \ref{connect}.

%% file: 7.0Acknowledgements.tex
\section*{Acknowledgement}

The authors would like to thank Guanghui Wang for his guidance and support throughout the preparation of this paper.


%% file: 8appendix.tex
\begin{appendix}

\section{Proof of Lemma \ref{findunit}}\label{appfindunit}
\begin{proof}[Proof of Lemma \ref{findunit}]
Let $G_0:=G\setminus W-U$. There exist $0<\epsilon_1 \leq 1/8,0<\epsilon_2<1$ and sufficiently large $C>0$ such that $m>8/\eta$ and Inequality (\ref{m200}) holds. Then 
$|W|\leq d^2 m^{y+1}/3\leq m^{y+1}dn/(3m^{200})\leq \eta(1-\eta)dn/2$. 

\begin{claim}
For any integers $s \leq t+r<200$, we can find vertex-disjoint $S(v_i)$-stars centered at $v_i$, $i= 1,\ldots, m^{s}$, each of size $(1-3\eta)d$ and $S(u_j)$-stars centered at $u_j$, $j= 1,\ldots, dm^t$, each of size $m^{r}$.
\end{claim} 

\begin{proof}
Indeed, let $S$ be the vertex set of a maximal collection $\mathcal{S}$ of vertex-disjoint stars constructed as above. If $\mathcal{S}$ is not as desired , then by Inequality (\ref{m200}) we have 
\[|S|\leq (d-3\eta d+1)m^s+(dm^{t}+1)m^r\leq 4 dm^{t+r}
   \leq 4 dm^{199}. \]
So we have \[|S\cup U|\leq \frac{dm^{y+1}}{4}+4dm^{199}\leq 5dm^{199}\leq \eta n.\]
By Proposition \ref{property}, we have \[d(G_0-S)=d(G\setminus W-(S\cup U)) \geq (1-3\eta)d.\] 
This implies that we can find one more star as desired, contradicting with the maximality of $\mathcal{S}$.
\end{proof}
Fix integers $s=y+3,t=y+5$ and $r=55$. 	
Let $V=\{v_1,v_2,\dots,v_{m^s}\}$ and $L(v_i)$ be the set of all leaves in each $S(v_i)$-star. Now, we shall use these vertex-disjoint stars to construct a $(d',m^y,m)$-unit in $G_0$ as follows.

\stepcounter{propcounter}
\begin{enumerate}[label=({\bfseries\Alph{propcounter}\arabic{enumi}})]
\rm\item\label{cunit1} Connect as many $(v_i,u_j)$ pairs as possible via the $(L(v_i),u_j)$-paths of length at most $m$ such that there is at most one path between any distinct pairs. 
\rm\item \label{cunit2} For each $v_i$, a leaf is \emph{occupied} if it is previously used as an endpoint of a path in \ref{cunit1}.			
\rm\item \label{cunit3} Let $(v_i,u_j)$ be the current pair to connect. Then the desired $(L(v_i),u_j)$-path shall avoid using
\begin{enumerate}[label=({\alph{propcounter}\arabic*})]
\rm\item \label{u1} any leaf occupied in $L(v_i)$ as an endpoint;
\rm\item \label{u2} edges in previous connecting paths, $\cup_{i=1}^{m^s} S(v_i)$ 
$\cup_{j=1}^{dm^t} S(u_j)$ and $W$;
\rm\item \label{u3} all vertices in $U\cup V$.
\end{enumerate}
\end{enumerate}

Firstly, we can get the following claim. 		
\begin{claim}\label{vicpuj}
There is a vertex $v_i$ connecting at least $q=d'+\eta d/2$ distinct centers $u_j$ satisfying \ref{cunit1}-\ref{cunit3}. 
\end{claim} 
\begin{proof}
Suppose to the contrary that each $v_i$ is connected to less than $q$ centers $u_j$. Then the number of vertices used in all $(L(v_i),u_j)$-paths for $i\in[m^s]$ is at most $q\cdot m\cdot m^{s}\leq dm^{s+1}$. There are at least $dm^t-dm^{s+1}\geq dm^t/2$ $S(u_j)$-stars that are completely disjoint from all those paths, which means that there are at least 
\[dm^t/2\geq dm^{y+2}:=x_0\]
available centers $u_j$, say $R'$. For each $S(v_i)$-star, there are at least $d-3\eta d-q=\eta d/2$ leaves not occupied. Thus, there is a set $L'$ of at least $\eta d/2\cdot m^s>x_0$ leaves not occupied from all $S(v_i)$-stars. Let $W'$ be the set of all the edges that we need to avoid in \ref{u2}. Then 
\[|W'|\leq dm^{s+1}+(d-3\eta d)m^s+dm^t \cdot m^r+|W|\leq \frac{d^2 m^{y+1}}{2}\leq \frac{d\rho(x_0)x_0}{2},\]
where the penultimate inequality holds because of Inequality (\ref{m60}).
Let $U'$ be the set of vertices that we need to avoid in \ref{u3}. Then
\[|U'|\leq |U|+|V|\leq\frac{ (m-1)dm^{y}}{4}+m^s\leq\frac{dm^{y+1}}{4}\leq \frac{\rho(x_0) x_0}{4}.\]
By Lemma \ref{robust}, we know that $G$ is an $(\epsilon_1,\epsilon_2 d)$-robust-expander. Applying Lemma \ref{rbsdl} to graph $G$ with $X_1=L'$, $X_2=R'$, $U=U'$ and $W=W'$, there is an $(L',R')$-path of length at most $m$, resulting in one more path between some $v_i, i\in[m^s]$ and $u_j, j\in[dm^t]$, a contradiction.
\end{proof}	

Without loss of generality, let $v_i,u_1,u_2,\dots,u_q$ be the centers guaranteed by Claim \ref{vicpuj}. 
Let $P$ be the interior vertex set of all paths between $v_i$ and $u_j$, $j\in[q]$. If $P$ is disjoint from $\cup_{j=1}^{q} S(u_j)$, then we get a desired unit.
Otherwise, we discard the $S(u_j)$-stars if at least half of its leaves are used in $P$. Since $|P|\leq qm$, the number of stars discarded is at most $\frac{qm}{m^r/2}\leq \eta d/2$. Hence, there are at least $d'$ $S(u_j)$-stars left. Observe that each $S(u_j)$-star left has at least $m^r/2\geq m^y$ leaves that do not appear in $P$. These $S(u_j)$-stars together with their corresponding paths to $v_i$ form a desired $(d',m^y,m)$-unit in $G_0$.
\end{proof}	
	
	\section{Proof of Lemma \ref{fu}}\label{appfu}
	\par The construction of units here essentially follows the arguments in Appendix \ref{appfindunit}.
	\begin{proof}[Proof of Lemma \ref{fu}]
		Let $G_1:=G\setminus W-U$ and $x,y$ be two positive integers such that  $x<100$ and $y=\lfloor (x-7)/2\rfloor$. 
		
		\begin{claim}
			There are vertex-disjoint stars in $G_1$, say $S(a_i)$ centered at $a_i$, $i\in\{1,2,\dots,m^{y+3}\}$, each of size $d/4$, and $S(b_j)$ centered at $b_j$, $j\in\{1,2,\dots,dm^{y+5}\}$, each of size $m^{y+1}$.
		\end{claim} 
		\begin{proof}
			Let $\mathcal{S}$ be the maximal collection of desired vertex-disjoint stars. Denote by $S$ the vertex set of $\mathcal{S}$. As $y=\lfloor (x-7)/2\rfloor$, if $\mathcal{S}$ is not as desired, then 
			\[|S|\leq (\tfrac{d}{4}+1)\cdot m^{y+3}+(m^{y+1}+1)\cdot dm^{y+5}\leq \frac{dm^x}{2}.\]
			Since $G$ is $(dm^x, d^2 m^x, d/2)$-dense and
			\[|S\cup U|\leq \frac{dm^x}{2}+dm^{y+1}\leq dm^x,\]
			we have $d(G\setminus W-(S\cup U))\geq d/2$. In addition, we know that $m^{y+1}\leq m^{46}\leq d/2$ using Inequality \ref{m60}. This implies that there exists one more star as desired, contradicting with the maximality of $\mathcal{S}$.
		\end{proof}
		
		Let $A=\{a_1,a_2,\dots,a_{m^{y+3}}\}$ and $L(a_i)$ be the set of all leaves in each star $S(a_i)$. Now, we shall use these vertex-disjoint stars to construct a $(3d/20,m^y,m)$-unit in $G_1$ as follows.
		\stepcounter{propcounter}
		\begin{enumerate}[label = ({\bfseries \Alph{propcounter}\arabic{enumi}})]
			\rm\item\label{edgevertex1} Connect as many $(a_i,b_j)$ pairs as possible via the $(L(a_i),b_j)$-paths of length less than $m$ such that there is at most one path between any pair.
			\rm\item\label{edgevertex2} For each $L(a_i)$, a leaf is \emph{occupied} if it is previously used as an endpoint of a path in \ref{edgevertex1}.
			\rm\item\label{edgevertex3} Let $(a_i,b_j)$ be the current pair to connect. Then an new  $(L(a_i),b_j)$-path shall avoid using
			\begin{itemize}
				\item[$\bullet$] any leaf occupied by $L(a_i)$ as an endpoint;
				\item[$\bullet$] edges in previous connecting paths, $\bigcup_{i=1}^{m^{y+3}} S(a_i)$, $\bigcup_{j=1}^{dm^{y+5}} S(b_j)$ and $W$;
				\item[$\bullet$] all vertices in $U\cup A$.
			\end{itemize}
		\end{enumerate}
		Then we have the following claim.
		
		\begin{claim}\label{123}
			There is a vertex $a_i$ connecting at least $q=d/5$ distinct centers $b_j$, and all these $(L(a_i),b_j)$-paths are pairwise edge disjoint.
		\end{claim} 
		
		\begin{proof}
			Suppose to the contrary that each $a_i$ is connected to less than $q$ centers $b_j$. Then the number of vertices used in all $(L(a_i), \{b_j\})$-paths is at most $q\cdot m\cdot m^{y+3}\leq dm^{y+4}$. Then there are at least $dm^{y+5}-dm^{y+4}\geq dm^{y+5}/2$ $S(b_j)$ stars that are completely disjoint from all those paths, which means that there are at least $dm^{y+5}/2\geq dm^{y+4}$ available centers $b_j$, say $R'$.  For each $S(a_i)$ star, there are at least $d/4-q\geq d/20$ leaves not occupied. Thus, there is a set $L'$ of at least $d/20\cdot m^{y+3}$ leaves not occupied by any $S(a_i)$ stars.  Set $x_0=d/20\cdot m^{y+3}$.
			
			Let $W'$ be the set of all the edges that we need to avoid in \ref{edgevertex3}. Then
			\[|W'|\leq dm^{y+4}+\tfrac{d}{4}\cdot m^{y+3}+dm^{y+5}\cdot m^{y+1}+|W|\leq d^2 m^{y+1}\leq d\rho(x_0)x_0.\]
			Let $U'$ be the set of vertices that we need to avoid in \ref{edgevertex3}. Then
			\[|U'|\leq |U|+|A|\leq dm^{y+1}+m^{y+3}\leq \frac{\rho(x_0) x_0}{4}.\]
			Applying Lemma \ref{rbsdl} to graph $G$ with $X_1=L'$, $X_2=R'$, $U=U'$ and $W=W'$, we obtain an $(L',R')$-path of length at most $m$, resulting in one more pair $a_i,b_j$ to be connected, a contradiction.
		\end{proof}
		
		Without loss of generality, let $a_i,b_1,b_2,\dots,b_q$ be the centers guaranteed by Claim \ref{123}. Let $P$ be the interior vertex set of all paths between $a_i$ and $b_j$, $j\in[q]$. If $P$ is disjoint from $\cup_{j=1}^{q} S(b_j)$, then we get a desired unit.
		Otherwise, we discard a $S(b_j)$-star if at least half of its leaves are used in $P$. Since $|P|\leq qm$, the number of stars discarded is at most $qm/(m^{y+1}/2)\leq d/20$. Hence, there are at least $q-d/20=3d/20$ $S(b_j)$-stars left. Note that each $S(b_j)$-star left has at least $m^{y+1}/2\geq m^y$ leaves that do not appear in $P$. These $S(b_j)$-stars together with their corresponding paths to $a_i$ form a $(3d/20,m^y,m)$-unit in $G_1$.
	\end{proof}
	\section{Proof of Lemma \ref{adjuster}}\label{appadjuster}
	\par
	To prove Lemma \ref{adjuster}, we need the following lemma to find some simple vertex-disjoint $(D,m,1)$-adjusters so as to link them up to obtain a desired $(D,m,r)$-adjuster.
	
	\begin{lemma}\label{1-adjuster}
		For each $0<\epsilon_1,\epsilon_2<1$, there exists $K>0$ such that the following holds for all $n$ and $d$ with $\log^{200}\leq d\leq n/K$. Let $x,y$ be two positive integers such that $1<y<x/2-3<28$ and $D=dm^y$. 
		If $G$ is an $n$-vertex $(dm^x, d^2 m^x, d/2)$-dense $(\epsilon_1,\epsilon_2 d)$-robust-expander graph with $\delta(G)=d$, then for any subset $U\subseteq V(G)$ with $|U|\leq 10D$, and any $W\subseteq E(G)$ with $|W|\leq d^2 m^{y-1}$, $G\setminus W-U$ contains a $(D, m/4, 1)$-adjuster.
	\end{lemma}
	
	\begin{proof}[Proof of Lemma \ref{adjuster}]
		Suppose $G\setminus W-U$ contains a $(D,m,r)$-adjuster for some maximal integer $1\leq r\leq dm^{y-3}$, say $\mathcal{A}_1:=(v_1,F_1,v_2,F_2,A_1)$. Note that such an adjuster exists from Lemma~\ref{1-adjuster}. Let $U_1=U\cup V(F_1)\cup V(F_2)\cup A_1$. Then $|U_1|\leq 4D$, and $|W|\leq d^2 m^{y-1}$. By Lemma \ref{1-adjuster}, there is a $(D,m/4,1)$-adjuster $\mathcal{A}_2:=(v_3,F_3,v_4,F_4,A_2)$ in $G\setminus W-U_1$. As $|F_1\cup F_2|=|F_3\cup F_4|=2D$, $|U\cup A_1\cup A_2|\leq D/(10m)+10mr+10m/4\leq dm^{y-1}/2\leq \rho(2D)2D/4$, and $|W|\leq d\rho(2D)2D$, by Lemma \ref{rbsdl}, there is a $(F_1\cup F_2, F_3\cup F_4)$-path $P'$ of length at most $m$, without loss of generality, we can say $P'$ is a $(F_1,F_3)$-path. Since  $F_1$ and $F_3$ are $(D,m)$-expansion of $v_1$ and $v_3$, respectively, then $P'$ can be extended to be a $(v_1,v_3)$-path $P$ of length at most $3m$. 
		
		We claim that $(v_2,F_2,v_4,F_4,A_1\cup A_2\cup P)$ is a $(D,m,r+1)$-adjuster. Indeed, we easily have that \ref{adjla1} and \ref{adjla2} hold, and $|A_1\cup A_2\cup P|\leq 10mr+10\cdot m/4+3m\leq 10m(r+1)$, so that \ref{adjla3} holds. Finally, let $\ell=\ell(\mathcal{A}_1)+\ell(\mathcal{A}_2)+\ell(P)$. If $i\in\{0,1,\ldots,r+1\}$, then there is some $i_1\in\{0,1,\ldots,r\}$ and $i_2\in\{0,1\}$ such that $i=i_1+i_2$. Let $P_1$ be a $(v_1,v_2)$-path of length $\ell(\mathcal{A}_1)+2i_1$ in $G\setminus W[A_1\cup\{v_1,v_2\}]$ and $P_2$ be a $(v_3,v_4)$-path of length $\ell(\mathcal{A}_2)+2i_2$ in $G\setminus W[A_2\cup\{v_3,v_4\}]$. Thus, $P_1\cup P\cup P_2$ is a $(v_2,v_4)$-path of length $\ell+2i$ in $G\setminus W[A_1\cup A_2\cup V(P)]$, and so $\ell$ satisfies \ref{adjla4}.
	\end{proof}
	
	Now, we only need to prove Lemma \ref{1-adjuster}. We use the idea of \cite{Luan2022BalancedSO}, and the following concept of octopus is also from \cite{Luan2022BalancedSO}.
	\begin{definition}
		\rm Given $r_1,r_2,r_3,r_4\in \mathbb{N}$, an \emph{$(r_1,r_2,r_3,r_4)$-octopus} $\mathcal{O}=(\mathcal{A},F,\mathcal{B},\mathcal{P})$ is a graph consisting of \begin{enumerate}
			\item[$\bullet$] a core $(r_1,r_2,1)$-adjuster $\mathcal{A}$, $F$ is one of the ends of $\mathcal{A}$, and
			\item[$\bullet$] a family of $r_3$ vertex-disjoint $(r_1,r_2,1)$-adjusters, denoted by $\mathcal{B}=\{\mathcal{A}_1,\ldots,\mathcal{A}_{r_3}\}$, which admits $\left(\bigcup_{i\in[r_3]}V(\mathcal{A}_i)\right)\cap V(\mathcal{A})=\emptyset$, and
			\item[$\bullet$] a minimal family $\mathcal{P}$ of internally vertex-disjoint paths of length at most $r_4$ such that at least one end of each adjuster in $\mathcal{B}$ is connected to $F$ by a subpath of $P$, where $P\in \mathcal{P}$, and all such subpaths are disjoint from all center sets of the adjusters in $\mathcal{B}\cup \mathcal{A}$.
		\end{enumerate}
	\end{definition}
	\begin{proof}[Proof of Lemma \ref{1-adjuster}]
		Let $G_1:=G\setminus W-U$. We claim that there are at least $m^x$ pairwise vertex-disjoint $(d/80,m/40,1)$-adjusters in $G_1$. Indeed, we may assume for contradiction that there are less than $m^{x}$ pairwise vertex-disjoint $(d/80,m/40,1)$-adjusters in $G_1$, and let $U_0'$ be the vertex set of all such adjusters. Then $|U_0'|\leq m^{x}(2\cdot d/80+10\cdot m/40)\leq dm^{x}/2$. Since $G$ is $(dm^x, d^2 m^x, d/2)$-dense,
		$|U\cup U_0'|\leq 10D+dm^{x}/2\leq dm^x$
		and $|W|\leq d^2 m^x$, then $d(G_1-U_0')\geq d/2$. By Corollary \ref{col}, there is a bipartite $(\epsilon_1,\epsilon_2 d)$-robust-expander $G'\subseteq G_1-U_0'$ with $\delta(G')\geq d/16$ for some $0<\epsilon_1,\epsilon_2<0$. Then by Bondy-Simonovits theorem, $ex(n,C_{2k})\leq O(n^{1+1/k})$, we have that the length of the shortest even cycle $C$ in $G'$ is at most $m/16$, and denote by $2r$ the length of $C$. Now, we arbitrarily choose $v_1,v_2\in V(C)$ of distance $r-1$ apart on $C$. Since $\delta(G')\geq d/16$, the order of $N_{G'-C}(v_1)$ and $N_{G'-C}(v_2)$ are both larger than $d/40$, where we also use Inequality (\ref{m60}). So we can choose $d/80$ distinct vertices from $N_{G'-C}(v_1)$ and $N_{G'-C}(v_2)$, respectively, together with $C$ forming a $(d/80,m/40,1)$-adjuster, a contradiction. 
		
Let $\mathcal{H}$ be a family of $m^x$ pairwise vertex-disjoint $(d/80,m/40,1)$-adjusters in $G_1$ and $U_0$ be the vertex set of all adjusters in $\mathcal{H}$. We say an adjuster is \emph{touched} by a path if they intersect in at least one vertex, and \emph{untouched} otherwise. Now we claim that there is a vertex set that can connect many ends of different adjusters through internally vertex-disjoint short paths.
		
		\begin{claim}\label{link}
			For any positive integer $t$ with $ t>y+1$, let $X\subseteq V(G)$ be an arbitrary set with $|X|\leq dm^{t}/2$, and $W$ be an edge set with $|W|\leq d^2 m^{t}$. Let $Y\subseteq V(G)-U$ with $|Y|\geq dm^{t+1}/80$, and $\mathcal{H}$ be a family of $(d/80,m/40,1)$-adjusters with $|\mathcal{H}|\geq m^{2t+1}$ in $G-(X\cup Y)$. Let $\mathcal{P}_{Y}$ be a maximum collection of internally vertex-disjoint paths of length at most $m/8$ in $G-X$, where each path connecting $Y$ to one end from distinct adjusters in $\mathcal{H}$. Then $Y$ can be connected to $1600m^{t+y+1}$ ends from distinct adjusters in $\mathcal{H}$ via some subpaths of the paths in $\mathcal{P}_{Y}$.
		\end{claim}
		\begin{proof}
			Suppose to the contrary that $Y$ is connected to less than $1600m^{t+y+1}$ ends from distinct adjusters, and denote by $P$ the set of internal vertices of those paths constructed as above. Then $|P|\leq 1600m^{t+y+1}\cdot m/8=200m^{t+y+2}$. 
			From the assumptions, there are at least $m^{2t+1}-200m^{t+y+2}=m^{t+1}(m^{t}-200m^{y+1})\geq m^{t+1}$ adjusters in $\mathcal{H}$ untouched by the paths in $\mathcal{P}_{Y}$.
			Arbitrarily choose  $m^{t+1}$ such adjusters, and let $E$ be the vertex set of the union of their ends. We get $|E|=2\cdot m^{t+1}\cdot d/80\geq dm^{t+1}/80=:D'$.
			Note that $|X\cup P|\leq dm^{t}/2+200m^{t+y+2}\leq dm^{t} \leq \rho(D') D'/4$, $|W|\leq d\rho(D') D'$ and $|Y|\geq D'$, there is a path of length at most $m/8$ between  $Y$ and $E$ by Lemma \ref{rbsdl}, a contradiction to the maximality of $\mathcal{P}_Y$.
		\end{proof}
		
		Now, we are ready to construct many octopuses via the above-mentioned small adjusters and paths.
		Let $B$ be the union of the center sets and core vertices of all those adjusters from $\mathcal{H}$.
		
		\begin{claim}\label{manyoctopus}
			Let $z$ be a positive integer such that $y+1<z<x/2-1<30$. Then $G_1$ contains $m^{z}$ $(d/80,m/40,800m^{y},m/8)$-octopus $\mathcal{O}_j=(\mathcal{A}_j,R_j,\mathcal{B}_j,\mathcal{P}_j)$, $1\leq j\leq m^{z}$ such that the following rules hold.
			\stepcounter{propcounter}
			\begin{enumerate}[label = ({\bfseries \Alph{propcounter}\arabic{enumi}})]
				\rm\item\label{ocolab1} For each $j\in[m^{z}]$, $\mathcal{A}_j$ are pairwise disjoint.
				\rm\item\label{ocolab2} $\mathcal{A}_i\notin \mathcal{B}_j$, $1\leq i,j\leq m^{z}$.
				\rm\item\label{ocolab3} Every adjuster in $\mathcal{B}_j$ intersects at least one path in $\mathcal{P}_j$, $1\leq j\leq m^{z}$.
				\rm\item\label{ocolab4} Each path $P$ in $\mathcal{P}_i$ satisfying $(B\cup V(\mathcal{A}_j))\cap P=\emptyset$, where $1\leq i\neq j\leq m^{z}$.
				\rm\item\label{ocolab5} For any $P\in \mathcal{P}_i$, $P'\in \mathcal{P}_j$, $V(P)\cap V(P')=\emptyset$, where $1\leq i< j\leq m^{z}$.
			\end{enumerate}
		\end{claim}
		
		\begin{proof}
			We construct the desired octopuses iteratively. Suppose that we can only construct less than $m^z$ octopuses. Let $U_1=U\cup B$. Then $|U_1|\leq |U|+|B|\leq 10D+(10\cdot m/40+2)m^x\leq 12D$. Let $U_0''$ be the union of the vertex set of the ends of the core adjusters of octopuses we have constructed. Then $|U_0''|< m^z\cdot 2\cdot d/80=dm^z/40$. We say an adjuster is \emph{used} if it appeared in previous octopuses, and \emph{unused} otherwise. Up to now, there are less than $m^z(800m^y+1)\leq 810m^{z+y}$ used adjusters, and so there are at least $m^x-810m^{y+z}\geq m^{2z+1}$ unused adjusters. Let $P^*$ be a set of all vertices in all $\mathcal{P}_j$, $j<m^z$. Then $|P^*|< m/8\cdot 800 m^y\cdot m^z=100 m^{z+y+1}$.
			
			Arbitrarily choose $m^{z+1}$ unused adjusters, denoted by $\mathcal{C}$. Let $X$ be the union of ends of the adjusters in $\mathcal{C}$. Then $|X|=m^{z+1}\cdot 2\cdot d/80\geq dm^{z+1}/80$. Note that there are at least $m^{2z+1}-m^{z+1}\geq m^{2z}$ unused adjusters remained apart from $\mathcal{C}$, denoted by $\mathcal{D}$. Let $U_2=U_0''\cup U_1\cup P^*$. Then we have $|U_2|\leq dm^z/40+12 D+100 m^{z+y+1}\leq dm^z/2$. Applying Claim \ref{link} with $(X,Y,\mathcal{H},t,W)=(U_2,X, \mathcal{D},z, W)$, we get that $X$ can connect to $1600m^{z+y+1}$ ends from different adjusters in $\mathcal{D}$ via some internally vertex-disjoint paths of length at most $m/8$. By the pigeonhole principle, there is an adjuster in $\mathcal{C}$, say $\mathcal{A}^*$, such that $\mathcal{A}^*$ has an end $R^*$ connected to at least $800 m^y$ adjusters, say $\mathcal{B}^*$, via a subfamily of internally vertex-disjoint paths, denoted by $\mathcal{P}^*$. From the construction process it is easy to see that \ref{ocolab1}-\ref{ocolab5} hold, and $\mathcal{A}^*$, $R^*$, $\mathcal{B}^*$ and $\mathcal{P}^*$ form one more $(d/80,m/40,800 m^{y},m/8)$-octopus.
		\end{proof}
		
		Now, we have $m^z$ $(d/80,m/40,800m^{y},m/8)$-octopus $\mathcal{O}_j=(\mathcal{A}_j,R_j,\mathcal{B}_j,\mathcal{P}_j)$, $1\leq j\leq m^{z}$. Let $L_j$ be the other end of $\mathcal{A}_j$ than $R_j$, and denote by $X'$ the union of all ends $L_j$, $1\leq j\leq m^z$. Then $|X'|=m^z\cdot d/80=dm^z/80$. Recall that there are $m^x$ adjusters and at most $m^z\cdot (800m^y+1)$ used adjusters. Thus, there are at least $m^x-m^z(800m^y+1)\geq m^{2z+1}$ unused adjusters, and denote by $\mathcal{D}'$ a collection of $m^{2z+1}$ unused adjusters among them. Let $P^{**}=\cup_{j=1}^{m^z} V(\mathcal{P}_j)$. Then $|P^{**}|\leq m/8\cdot 800 m^y\cdot m^z\leq 100m^{z+y+1}$. By definition, for each octopus $\mathcal{O}_j=(\mathcal{A}_j,R_j,\mathcal{B}_j,\mathcal{P}_j)$, $1\leq j\leq m^{z}$, every adjuster $\mathcal{A}\in \mathcal{B}_j$ intersects $V(\mathcal{P}_j)$ and thus there is a shortest path in $\mathcal{A}$ of length at most 2 connecting a core vertex of $\mathcal{A}$ to $V(\mathcal{P}_j)$, and denote by $\mathcal{Q}_j$ the disjoint union of such paths take over all adjusters in $\mathcal{B}_j$. Let $Q'=\cup_{j=1}^{m^z}V(\mathcal{Q}_j)$. Then $|Q'|\leq 3\cdot |\mathcal{B}_j|\cdot m^z\leq m^{z+y+1}$. Let $U_2'=U\cup B\cup P^{**}\cup Q'$. Then $|U_2'|\leq 12 D+100 m^{z+y+1}+m^{z+y+1}\leq dm^{z}/2$. By the same approach, applying Claim \ref{link} with $(X,Y,\mathcal{H},t,W)=(U_2',X', \mathcal{D}',z, W)$, we know that $X$ can be connected to $1600 m^{z+y+1}$ ends from different adjusters in $\mathcal{D}'$ via some internally vertex-disjoint paths of length at most $m/8$. By pigeonhole principle, there is an core adjuster $\mathcal{A}_k$ such that $L_{k}$ is connected to a family $\mathcal{B}_k'$ of at least $800 m^y$ adjusters, via a subfamily of internally vertex-disjoint paths, denoted by $\mathcal{P}_{k}'$. Then  $\mathcal{A}_{k}$, $L_{k}$, $\mathcal{B}_{k}'$ and $\mathcal{P}_{k}'$ form a $(d/80,m/40,800 m^{y},m/8)$-octopus. Note that $\mathcal{A}_{k}$, $R_{k}$, $\mathcal{B}_{k}$ and $\mathcal{P}_{k}$ also form a $(d/80,m/40, 800m^{y},m/8)$-octopus. 
		
		Denote by $A_k$ the center vertex set of the adjuster $\mathcal{A}_k$. Recall that $\mathcal{A}_k$ is a $(d/80,m/40,1)$-adjuster, so $L_k$ and $R_k$ are $(d/80,m/40)$-expansions of vertices $v_1, v_2$ respectively. Let $F_1'= G\setminus W[V(L_k)\cup V(\mathcal{B}_k')\cup V(\mathcal{P}_k')]$ and $F_2'$ is a component of $G\setminus W[V(R_k)\cup V(\mathcal{B}_k)\cup V(\mathcal{P}_k)]-V(\mathcal{P}_k')$ containing $v_2$. Note that all paths in $\mathcal{P}_k$ and $\mathcal{P}_k'$ are disjoint from $B$, and $V(\mathcal{P}_k)$ and $V(\mathcal{P}_k')$ are disjoint. By the fact that $V(\mathcal{P}_k')$ is disjoint from $B$ and $Q$, we get 
		$|F_2'|\geq |V(\mathcal{B}_k)|-|V(\mathcal{P}_k')|\geq 800m^y\cdot 2\cdot d/80-800m^y\cdot m/8 \geq 3D$, and the distance between $v_2$ and each $v\in V(F_2')$ is at most $m/40+m/8+m/40+m/32+m/40\leq m/4$. By Proposition \ref{DD}, there is a subgraph $F_2\subseteq F_2'$, which is a $(D,m/4)$-expansion of $v_2$. Similarly, we can find a subgraph $F_1\subseteq F_1'$, which is a $(D,m/4)$-expansion of $v_1$, and $F_1$ and $F_2$ are disjoint. Recall that $C\cup \{v_1, v_2\}$ is an even cycle of length $2r\leq m/16$ and the distance between $v_1$ and $v_2$ on $C\cup\{v_1,v_2\}$ is $r-1$. Hence, $(v_1,F_1,v_2,F_2,C)$ is a $(D,m/4,1)$-adjuster.
		
	\end{proof}
	\section{Proof of Lemma \ref{fp}}\label{appfp}
	\begin{proof}[Proof of Lemma \ref{fp}]
		Let $G_1:=G\setminus W-U$.
		Note that $|U|\leq dm^{y-2}\leq \rho(D)D/4$, $|W|\leq d^2m^{y-2}\leq d\rho(D)D$, $|Z_1\cup Z_2|\geq D$ and $|I_1\cup I_2|\geq D$. By Lemma \ref{rbsdl}, there is a $(Z_1\cup Z_2, I_1\cup I_2)$-path $P_1$ of length at most $m$ in $G_1$, say $P_1\cap Z_1={z_1}$. Since $I_1$ is a $(D,m)$-expansion of $v_1$, $P_1$ can be extended to a $(z_1,v_1)$-path $P$ of length at most $2m$ where $z_1\in Z_1$.
		
		Now, we need to find a $(z_2,v_2)$-path $Q$ in $G_1-P$ such that $\ell\leq\ell(P)+\ell(Q)\leq\ell+10m$ where  $z_2\in Z_2$ and $z_1\neq z_2$. Let $(X,v',I')$ be a triple set such that $\ell(X)$ is maximised and satisfies the following properties.
		\begin{itemize}
			\item[(P1)] $I'$ is a $(D,m)$-expansion of $v'$ in $G_1-P$.
			\item[(P2)] $X$ is a $(v',v_2)$-path in $G_1-P$ and $I'\cap X={v'}$.
			\item[(P3)] $\ell(X)\leq\ell+4m.$
		\end{itemize}
		Note that such triple set must exist because of the basic case where $I'=I_2, v'=v_2$ and $X=G[v']$. 
		
		We claim that $\ell(X)\geq \ell$. Otherwise, we can get a contradiction. Let $U_1:=V(P)\cup V(X)\cup V(I_2)\cup V(I')$. Then 
		\[|U\cup U_1|\leq dm^{y-2}+2m+\ell +D+D\leq 5D\leq dm^{y+1}.\]
		Applying Lemma \ref{fu} to $G$ with $W=W$ and $U=U\cup U_1$, there is a $(3d/20,m^y,m)$-unit $Z_0$ with core vertex $z_0$ in $G_1-U_1$. Since $|U\cup V(X)\cup V(P)|\leq dm^{y-2}+\ell +2m\leq \rho(D)D/4$, applying Lemma \ref{rbsdl} to $G$ with $(U,W,X_1,X_2)=(U\cup X\cup P, W, Z_0, I')$, there is a $(Z_0,I')$-path $Y'$ of length at most $m$ in $G_1-X-P$. Recall that $Z_0$ is a $(3d/20,m^y,m)$-unit and $I'$ is a $(D,m)$-expansion, so $Y'$ can be extended to a $(z_0,v')$-path $Y$ of length at most $m+m+1+m\leq 4m$. By the property of $(3d/20,m^y,m)$-unit, we can find a subgraph $Z_0'\subseteq Z_0-Y+z_0$ which is a $(D,m)$-expansion of $z_0$. This implies there is a $(v_2,z_0)$-path $v_2Xv'Yz_0$, briefly denoted as $X'$. Since
		\[\ell(X)+1\leq \ell(X')\leq \ell(X)+4m,\]
		there is a triple set $(X',z_0,Z_0')$ satisfying the properties (P1)-(P3) with $\ell(X')>\ell(X)$, a contradiction to the maximality of $\ell(X)$.
		
		Applying Lemma \ref{rbsdl} to $G$ with $(U,W,X_1,X_2)=(U\cup X\cup P, W, Z_2, I')$, there is a $(Z_2,I')$-path $Q'$ of length at most $m$ in $G_1-X-P$. $Q'$ can be extended to a $(z_2,v')$-path $Q''$ with length at most $4m$. Let $Q:=z_2 Q'' v' X v_2$. Then $\ell\geq\ell(X)\leq\ell(Q)\leq \ell(X)+\ell(Q'')\leq \ell+4m+4m\leq \ell+8m$. 
		
		Thus we find two vertex-disjoint paths $P$ and $Q$ with $\ell\leq\ell(P)+\ell(Q)\leq \ell+8m+2m\leq\ell+10m$ such that $P$,$Q$ link $\{z_1,z_2\}$ to $\{v_1,v_2\}$ for $z_i\in Z_i$ and $z_1\neq z_2$.
	\end{proof}
\end{appendix}